\UseRawInputEncoding
\documentclass[12pt]{article}

\usepackage{mathrsfs}
\usepackage{amssymb}
\usepackage{amsmath}
\usepackage{amsfonts}
\usepackage{amstext}
\usepackage{amsthm}
\usepackage{graphicx}
\usepackage{subfig}
\usepackage{color}
\usepackage{indentfirst,latexsym,bm}
\usepackage{mathrsfs}
\usepackage{cite}
\usepackage[all]{xy}
\usepackage{newlfont}

\usepackage{enumerate}
\numberwithin{equation}{section}

\numberwithin{equation}{section}

\begin{document}

\title{Sign-changing solution for an overdetermined elliptic problem on unbounded domain
\thanks{Research was supported by the National Natural Science Foundation of China (Grants No. 12371110 and
No. 12301133).}}

\author{Guowei Dai\thanks{Corresponding author.
\newline
School of Mathematical Sciences, Dalian University of Technology, Dalian, 116024, P.R. China
\newline
\text{~~~~ E-mail}: daiguowei@dlut.edu.cn.}, Yong Zhang\thanks{Institute of Applied System Analysis, Jiangsu University, Zhenjiang, 212013, P.R. China
\newline
\text{~~~~ E-mail}: zhangyong@ujs.edu.cn} \\
}
\date{}
\maketitle

\renewcommand{\abstractname}{Abstract}

\begin{abstract}
We prove the existence of two smooth families of unbounded domains in $\mathbb{R}^{N+1}$ with $N\geq1$
such that
\begin{equation}
-\Delta u=\lambda u\,\, \text{in}\,\,\Omega, \,\, u=0,\,\,\partial_\nu u=\text{const}\,\,\text{on}\,\,\partial\Omega\nonumber
\end{equation}
admits a sign-changing solution.
The domains bifurcate from
the straight cylinder $B_1\times \mathbb{R}$, where $B_1$ is the unit ball in $\mathbb{R}^N$. These results can be regarded as counterexamples to the Berenstein conjecture on unbounded domain. Unlike most previous papers in this direction, a very delicate issue here is that there may be two-dimensional kernel space at some bifurcation point. Thus a Crandall-Rabinowitz type bifurcation theorem from high-dimensional kernel space is also established to achieve the goal.
\end{abstract}


\emph{AMS Subjection Classification(2020):} 35B05; 35B32; 35N05; 47J15

\tableofcontents

\section{Introduction}
\quad\, Overdetermined elliptic problem is related to plasma physics \cite{Temam},
nuclear reactors \cite{Berenstein} and tomography \cite{Shepp, Smith}. Since Serrin's famous work \cite{Serrin}, overdetermined elliptic problem has attracted the attention of many mathematicians.
In the past half century, many celebrated results have been established, for instance \cite{Aftalion, BCN, ChenLi, Gidas, Pucci, Reichel, RRS}. 
In particular, Berenstein \cite{Berenstein, Berenstein1, Berenstein2} proposed the following conjecture.
\\ \\
\textbf{Berenstein conjecture.}
\emph{Let $\Omega$ be a bounded $C^{2,\alpha}$ domain in $\mathbb{R}^{N+1}$ with $\alpha\in(0,1)$. If there exists a nontrivial
solution $u$ of the overdetermined eigenvalue problem
\begin{equation}\label{overdeterminedeigenvalueproblem}
\left\{
\begin{array}{ll}
\Delta u+\lambda u=0\,\, &\text{in}\,\, \Omega,\\
u=0 &\text{on}\,\, \partial \Omega,\\
\partial_\nu u=\text{const} &\text{on}\,\, \partial \Omega,
\end{array}
\right.
\end{equation}
then $\Omega$ is a ball, where $\nu$ is the unit outer normal vector on $\partial \Omega$.}\\

For $N=1$, Berenstein \cite{Berenstein} proved that the existence of infinitely many eigenvalues for (\ref{overdeterminedeigenvalueproblem}) is equivalent to $\Omega$ being a disk, which was extended
to the Poincare disk with the hyperbolic metric by Berenstein and Yang \cite{Berenstein1}. Then Berenstein and Yang \cite{Berenstein2} further showed that Berenstein's result is valid for any $N$.
In \cite{Liu}, Liu proved that the Berenstein
conjecture holds if and only if the second order interior normal derivative of the corresponding
Dirichlet eigenfunction $u$ is constant on the boundary of $\Omega$.

The other long standing open problem considered in \cite{Berenstein1,Berenstein2}, which is very similar to (\ref{overdeterminedeigenvalueproblem}), is the following Schiffer Conjecture.
\\ \\
\textbf{Schiffer Conjecture.} \emph{Let $\Omega\subset \mathbb{R}^N$ be a bounded regular domain. Assume $u: \Omega \rightarrow \mathbb{R}$
is a solution to the problem:
\begin{equation}\label{SchifferConjecture}
\left\{
\begin{array}{ll}
\Delta u+\lambda u=0\,\, &\text{in}\,\, \Omega,\\
u=1 &\text{on}\,\, \partial \Omega,\\
\partial_\nu u=0 &\text{on}\,\, \partial \Omega.
\end{array}
\right.
\end{equation}
Then $\Omega$ is a ball and $u$ is radially symmetric.}
\\

The Schiffer conjecture is closely connected with Pompeiu problem \cite{Pompeiu, Pompeiu1}.
Let $\Omega$ be a nonempty bounded open subset of $\mathbb{R}^N$ with $N\geq  2$ and let $\mathcal{M}$ denote the set of rigid
motions of $\mathbb{R}^N$ onto itself. The domain $\Omega$ is said to have the Pompeiu property if and only if $f\equiv 0$ on $\mathbb{R}^N$ is the only continuous function satisfying
\begin{equation}
\int_{\sigma(\Omega)}f(x)\,\text{d}x=0\nonumber
\end{equation}
for every $\sigma\in \mathcal{M}$. The Pompeiu problem originated from harmonic analysis, which consists in characterizing the class of domains in $\mathbb{R}^N$ with the Pompeiu property. It has been proved \cite{Williams} that the problem (\ref{SchifferConjecture}) admits a solution with $\lambda>0$ if and only if the smooth $\Omega$ with $\partial\Omega$ connected fails to possess the Pompeiu property. The Schiffer conjecture has also been included in Yau's famous list of problems \cite[Problem 80]{Yau}. So far, there are few results about this conjecture. In \cite{Berenstein1,Berenstein2}, Berenstein and Yang found that the existence of infinitely many eigenvalues to (\ref{SchifferConjecture}) implies that $\Omega$ must be a round ball. Recently, Fall, Minlend and Weth \cite{FallMW} constructed a nontrivial family of compact subdomains of the flat cylinder $\mathbb{R}^N\times \mathbb{R}/2\pi\mathbb{Z}$ such that (\ref{SchifferConjecture}) admits eigenfunctions. That is to say, they gave the first counterexample to the Schiffer conjecture on unbounded domains in some sense.

On the other hand, if we only take the positive solutions of (\ref{overdeterminedeigenvalueproblem}) into consideration, which is strongly related with Berestycki-Caffarelli-Nirenberg (BCN) conjecture\cite{BCN}. It is well known that many positive answers to the BCN conjecture have been given in bounded domain by using the so-called moving plane method. Assume that $\Omega$ is $C^3$ and uniformly Lipschitz epigraph of $\mathbb{R}^2$ or $\mathbb{R}^3$, Farina and Valdinoci \cite{Farina} proved that
there exists no solution $u\in C^2\left(\overline{\Omega}\right) \cap L^\infty (\Omega)$ of
\begin{equation}\label{positivesoution}
\left\{
\begin{array}{ll}
\Delta u+\lambda u=0\,\, &\text{in}\,\, \Omega,\\
u>0 &\text{in}\,\, \Omega,\\
u=0 &\text{on}\,\, \partial \Omega,\\
\partial_\nu u=\text{const} &\text{on}\,\, \partial \Omega.
\end{array}
\right.
\end{equation}
That is to say if (\ref{positivesoution}) has a solution, $\Omega$ must be half-space, which can be seen as a confirm answer to BCN or Berenstein conjecture on unbounded domain. The first counterexample to BCN conjecture on unbounded domain was constructed by Sicbaldi \cite{Sicbaldi} via showing that the cylinder $B_1\times \mathbb{R}$ with $N\geq2$ can be perturbed to an unbounded domain whose boundary is a periodic
hypersurface of revolution with respect to the $\mathbb{R}$-axis and such that problem (\ref{positivesoution}) has a bounded solution, where $B_1$ is the unit ball of $\mathbb{R}^N$ centered on the origin.
Subsequently, Schlenk and Sicbaldi \cite{Schlenk} further proved that the above conclusion is valid for $N=1$ and these new extremal domains belong to a smooth bifurcation family of domains. Nontrivial domain emanating from the half-space has been obtained in \cite{Del} in $N\geq8$, and bifurcating from the complement of a ball has been given in \cite{Ros1}.

These results on positive solutions have enriched counterexamples to the BCN conjecture on unbounded domain. The natural question is whether there is a nontrivial unbounded region on which (\ref{overdeterminedeigenvalueproblem}) has nonsymmetric sign-changing solution, which means that we can construct a counterexample to the Berenstein conjecture on an unbounded region. The main contribution of this paper is to give a confirm answer to the question as follows.
\\ \\
\textbf{Theorem 1.1.} \emph{Let} $\mathcal{C}^{2,\alpha}_{\text{even},0}\left(\mathbb{R}/2\pi \mathbb{Z}\right)$
\emph{be the space of even $2\pi$-periodic $\mathcal{C}^{2,\alpha}$ functions of mean zero.
For each $N\geq2$ there exists a positive number $T_*$ with
\begin{equation}
T_*\in\left(\frac{2\pi}{\sqrt{\lambda_2}},\frac{2\pi}{\sqrt{\lambda_2-\lambda_1}}\right),\nonumber
\end{equation}
where $\lambda_1$ and $\lambda_2$ are the first and second eigenvalues of the zero-Dirichlet Laplacian on the unit ball, and a smooth map}
\begin{equation}
(-\varepsilon,\varepsilon)\rightarrow \mathcal{C}^{2,\alpha}_{\text{even},0}\left(\mathbb{R}/2\pi \mathbb{Z}\right)\times \mathbb{R}:s\mapsto \left(w_s,T_s\right)\nonumber
\end{equation}
\emph{with $w_0 = 0$, $T_0 = T_*$ and such that for each $s\in  (-\varepsilon,\varepsilon)$, problem (\ref{overdeterminedeigenvalueproblem}) has a sign-changing $T_s$-periodic solution
$u_s \in  \mathcal{C}^{2,\alpha}\left(\Omega_s\right)$ on the modified cylinder
\begin{equation}
\Omega_s=\left\{(x,t)\in \mathbb{R}^N\times \mathbb{R}:\vert x\vert<1+s\cos \left(\frac{2\pi}{T_s}t\right)+s w_s\left(\frac{2\pi}{T_s}t\right)\right\}.\nonumber
\end{equation}
Moreover, there exists a unique $T_s$-periodic function $r\in  \mathcal{C}^{2,\alpha}\left(\mathbb{R}\right)$ such that
\begin{equation}
r(t):\mathbb{R}\longrightarrow \left(0,1+s\cos \left(\frac{2\pi}{T_s}t\right)+s w_s\left(\frac{2\pi}{T_s}t\right)\right)\nonumber
\end{equation}
and $u_s(r(t),t)=0$.}\\

From Theorem 1.1, we see that the obtained solution is sign-changing and the bifurcation point $T_*$ is bounded from the above. In addition, we also find another bifurcation point $T^*$ which may be corresponding to high-dimensional kernel space.
\\ \\
\textbf{Theorem 1.2.} \emph{Under the assumptions of Theorem 1.1, there exist a positive number $T^*$ with
\begin{equation}
T^*>\frac{2\pi}{\sqrt{\lambda_2-\lambda_1}}\nonumber
\end{equation}
and a smooth map}
\begin{equation}
(-\varepsilon,\varepsilon)\rightarrow \mathcal{C}^{2,\alpha}_{\text{even},0}\left(\mathbb{R}/2\pi \mathbb{Z}\right)\times \mathbb{R}:s\mapsto \left(w_s,T_s\right)\nonumber
\end{equation}
\emph{with $w_0 = 0$, $T_0 = T^*$ and such that for each $s\in  (-\varepsilon,\varepsilon)$, problem (\ref{overdeterminedeigenvalueproblem}) has a sign-changing $T_s$-periodic solution
$u_s \in  \mathcal{C}^{2,\alpha}\left(\Omega_s\right)$ on the modified cylinder
\begin{equation}
\Omega_s=\left\{(x,t)\in \mathbb{R}^N\times \mathbb{R}:\vert x\vert<1+s\left(\beta\cos \left(\frac{2\pi}{T_s}t\right)+\gamma\cos \left(\frac{2m\pi}{T_s}t\right)\right)+s w_s\left(\frac{2\pi}{T_s}t\right)\right\},\nonumber
\end{equation}
where either $\beta$, $\gamma$ are any given nonzero constants with $\beta^2+\gamma^2=1$ if there exists some $m\in \mathbb{N}$ such that $T^*=m T_*$, or $\beta=1$ and $\gamma=0$ if $T^*\neq m T_*$ for any $m\in \mathbb{N}$.
Moreover, there exists a unique $C^{2,\alpha}$ function
\begin{equation}
r(t):\mathbb{R}\longrightarrow \left(0,1+s\left(\beta\cos \left(\frac{2\pi}{T_s}t\right)+\gamma\cos \left(\frac{2m\pi}{T_s}t\right)\right)+s w_s\left(\frac{2\pi}{T_s}t\right)\right)
\nonumber\end{equation}
 such that
$u_s(r(t),t)=0$ and it is $T_s$-periodic.}
\\ \\
\textbf{Remark 1.3.}
\emph{For $N=1$, in Section 5, we will find that $T_*=4/3$ and $T^*=4\sqrt{5}/5$ by an exact computation, which implies that $\beta=1$ and $\gamma=0$ in this case.    However, for $N\geq2$, $T^*$ may be the integer multiple of $T_*$, which gives that $\beta$, $\gamma$ are nonzero constants with $\beta^2+\gamma^2=1$. That is to say, the corresponding kernel space is exactly two-dimensional,
which leads to that the classical Crandall-Rabinowitz bifurcation theorem cannot be used directly. Therefore, we also establish a new Crandall-Rabinowitz type bifurcation theorem with high-dimensional kernel space in Section 4.}\\

As shown in Theorem 1.1 and Theorem 1.2, the main goal of the present paper
is to construct a nontrivial unbounded domain such that problem (\ref{overdeterminedeigenvalueproblem}) admits a sign-changing solution. We choose the eigenfunction corresponding to the second eigenvalue of zero-Dirichlet Laplacian on the unit ball as the trivial sign-changing solution. In the spirit of local bifurcation, we successfully achieve this aim. As far as we know, this is a new way to construct the sign-changing solution to the overdetermined problems. This method may be applied in other settings, such as other operators, domains in Riemannian Manifolds. However, there may be some new difficulties that need to be overcome with more effort.
Here the Bessel functions play a crucial role. For other operators or domains in Riemannian Manifolds, Bessel functions may not be used directly.
Perhaps studying other hypergeometric functions (such as the Legendre functions) can work for investigating domains in Riemannian Manifolds.
By establishing the relationship between eigenvalues and the results of this paper, the expected properties of other operators may be derived. These questions are very interesting to be considered in future.

We also like to mention the recent work \cite{Minlend,Ruiz} due to Minlend and Ruiz respectively, where sign-changing solutions to some overdetermined problems are obtained by using the local bifurcation theorem. The choice in \cite{Ruiz} of trivial sign-changing solutions depend closely on the form of the equation in which the nonlinearity is the Allen-cahn type, which is indeed a delicate issue. Moreover, the nontrivial domains constructed in \cite{Ruiz} are bounded, which is fundamentally different from our results. Although the unbounded domains constructed in \cite{Minlend} are periodic in the first coordinate and they bifurcate from suitable strips in $\mathbb{R}^2$, the Neumann boundary condition considered in \cite{Minlend} is varying from top to bottom which is different from the Berenstein problem (\ref{overdeterminedeigenvalueproblem}).
Concretely, assume that $\Omega_s\subset \mathbb{R}^2$ is the domains obtained in \cite[Theorem 1.1]{Minlend} and
let $\partial\Omega_s^+$ and $\partial\Omega_s^-$ denote the top and bottom boundary of $\Omega_s$, respectively, i.e.,
\begin{equation}
\partial\Omega_s^+=\left\{(x,t)\in\partial\Omega_s:t>0\right\}\,\,\,\text{and}\,\,\,\partial\Omega_s^-=\left\{(x,t)\in\partial\Omega_s:t<0\right\}.\nonumber
\end{equation}
The solution $u$ obtained in \cite[Theorem 1.1]{Minlend} satisfies $\partial_\nu u=-1$ on $\partial\Omega_s^+$ and $\partial_\nu u=+1$ on $\partial\Omega_s^-$ which is different from ours because here the Neumann boundary condition is the same on the whole boundary of $\Omega_s$.
Anyway, these work on changing-sign solutions are important in application and will be the object of intensive research.

The rest of this paper is arranged as follows. In Section 2, we give some preliminaries. In Section 3, we study the properties of a certain eigenvalue, which is key to obtain Theorems 1.1--1.2. The Section 4 is devoted to completing the proofs of Theorems 1.1--1.2.
In the Last Section, we show that the conclusions of Theorems 1.1--1.2 are also valid for $N=1$, where the method used is different from the case of $N\geq2$. In an appendix, we give a supplementary proof of \cite[Claim 6.6]{Schlenk} by filling a small gap.

\section{Preliminaries}

\quad\,  The main strategy of this paper is to apply Crandall-Rabinowitz bifurcation theorem.
The first step is to transform the problem (\ref{overdeterminedeigenvalueproblem}) into an abstract operator equation.
To achieve this goal, we show two elementary results on zero-Dirichlet Laplacian eigenvalue problems on the unit ball and the cylinder in the following.

\subsection{Some results on eigenvalue problems}

\quad\, We first consider the following eigenvalue problem
\begin{equation}\label{eigenvalueonball1}
\left\{
\begin{array}{ll}
\Delta u+\lambda u=0\,\, &\text{in}\,\, B_1,\\
u=0 &\text{on}\,\, \partial B_1.
\end{array}
\right.
\end{equation}
It is well known (see \cite{Coddington, Ince} or
\cite[p. 269]{Walter}) that problem (\ref{eigenvalueonball1}) possesses a sequence eigenvalues $0<\lambda_1<\lambda_2<\cdots<\lambda_k\nearrow+\infty$ for $k\in \mathbb{N}$.
Let $\overline{\phi}_k$ be the corresponding radial eigenfunction
to $\lambda_k$ with $\int_{B_1}\overline{\phi}_k^2\,\text{d}x=1/\left(2\pi\right)$ and $\overline{\phi}_k(0)>0$.
From now on, we use $C$ or $C_k$ to denote positive constant whose exact value may change from line to line.

In particular, when $N=1$, we have that
\begin{equation}
\lambda_k=\frac{(2k-1)^2\pi^2}{4}\,\,\,\text{and}\,\,\,\overline{\phi}_k(r)=\frac{1}{\sqrt{2\pi}}\cos\left(\frac{(2k-1)\pi}{2}r\right).\nonumber
\end{equation}
While, for $N=3$, we have the following lemma.
\\ \\
\textbf{Lemma 2.1.} \emph{When $N=3$ and $k\in \mathbb{N}$, one has that}
\begin{equation}
\lambda_k=k^2\pi^2\,\,\,\text{and}\,\,\,C_k\overline{\phi}_k(r)=\left\{
\begin{array}{ll}
\frac{\sin(k\pi r)}{r}\,\, &\text{if}\,\, r>0,\\
k\pi &\text{if}\,\, r=0,
\end{array}
\right.\nonumber
\end{equation}
\emph{where $C_k$ is chosen such that $\int_{B_1}\overline{\phi}_k^2\,\text{d}x=1/\left(2\pi\right)$}.
\\ \\
\textbf{Proof.} Let
\begin{equation}
f(r)=\left\{
\begin{array}{ll}
\frac{\sin(\sqrt{\lambda} r)}{r}\,\, &\text{if}\,\, r>0,\\
\sqrt{\lambda} &\text{if}\,\, r=0.
\end{array}
\right.\nonumber
\end{equation}
Then, for $r>0$,  we have that
\begin{equation}
f'(r)=\sqrt{\lambda}\frac{\cos (\sqrt{\lambda} r)}{r}-\frac{\sin (\sqrt{\lambda} r)}{r^2}\nonumber
\end{equation}
and
\begin{equation}
f''(r)=\lambda\frac{-\sin (\sqrt{\lambda} r)}{r}-2\sqrt{\lambda}\frac{\cos (\sqrt{\lambda} r)}{r^2}+2\frac{\sin (\sqrt{\lambda} r)}{r^3}.\nonumber
\end{equation}
Hence we have that
\begin{eqnarray}
f''(r)+\frac{2}{r}f'(r)+\lambda f(r)&=&\lambda\frac{-\sin (\sqrt{\lambda} r)}{r}-2\sqrt{\lambda}\frac{\cos (\sqrt{\lambda} r)}{r^2}+2\frac{\sin (\sqrt{\lambda} r)}{r^3}\nonumber\\
& &+2\sqrt{\lambda}\frac{\cos (\sqrt{\lambda} r)}{r^2}-2\frac{\sin (\sqrt{\lambda} r)}{r^3}+\lambda\frac{\sin(\sqrt{\lambda} r)}{r}\nonumber\\
&=&0.\nonumber
\end{eqnarray}
Further, we can verify that $f'(0)=0$.
Combining this with the boundary condition, we obtain that
\begin{equation}
\lambda_k=k^2\pi^2\,\,\,\text{and}\,\,\,C_k\overline{\phi}_k(r)=\left\{
\begin{array}{ll}
\frac{\sin(k\pi r)}{r}\,\, &\text{if}\,\, r>0,\\
k\pi &\text{if}\,\, r=0.
\end{array}
\right.\nonumber\nonumber
\end{equation}
for $r\in[0,1]$.\qed\\

It is well known that, for the one-dimensional case, the eigenvalues and eigenfunctions can be calculated by explicit expressions, which indicates that the Weyl's asymptotic formula holds \cite{Chavel}.
The conclusion of Lemma 2.1 shows that, for the three-dimensional case, the eigenvalues and eigenfunctions can still be calculated by explicit expressions, which plays an important role in our subsequent argument. This result also indicates that the equality holds in the Weyl's asymptotic formula (in fact, the radial situation can be understood as one-dimensional problem) to
three-dimensional radial eigenvalue problem.

Note that $\left(\lambda_k, {\phi}_k(x,t)\right)$ with ${\phi}_k(x,t)=\overline{\phi}_k(x)$ is also a solution pair of
\begin{equation}\label{eigenvalueoncylinder1}
\left\{
\begin{array}{ll}
\Delta_{\mathring{g}} {\phi}+\lambda {\phi}=0\,\, &\text{in}\,\, C_1^T,\\
{\phi}=0 &\text{on}\,\, \partial C_1^T,
\end{array}
\right.
\end{equation}
where $\mathring{g}$ is the Euclidean metric and
\begin{equation}
C_{1}^T=\left\{(x,t)\in \mathbb{R}^N\times \mathbb{R}/T\mathbb{Z}: \vert x\vert<1\right\}.\nonumber
\end{equation}
That is to say, the eigenvalues and eigenfunctions of (\ref{eigenvalueonball1}) are also the eigenvalues and eigenfunctions of (\ref{eigenvalueoncylinder1}).
Clearly, we see that
\begin{equation}
\int_{C_1^{2 \pi}}{\phi}_k^2\,\text{d} \text{vol}_{\mathring{g}}=1.\nonumber
\end{equation}
Since $\phi_k$ does not depend on $t$ and is radial, we will denote $\phi_k(x,t)$ by $\phi_k(r)$ with $r=\vert x\vert$.

\subsection{Rephrasing the problem}

\quad\,We now transform the aim problem into some abstract operator equation.
For each $v\in \mathcal{C}^{2,\alpha}_{\text{even},0}\left(\mathbb{R}/2\pi \mathbb{Z}\right)$ with $v>-1$, define
\begin{equation}
C_{1+v}^T=\left\{(x,t)\in \mathbb{R}^N\times \mathbb{R}/T\mathbb{Z}: \vert x\vert<1+v\left(\frac{2\pi t}{T}\right)\right\}\nonumber
\end{equation}
for all $T>0$.
We consider the following eigenvalue problem
\begin{equation}\label{eigenvalueproblem}
\left\{
\begin{array}{ll}
\Delta_{\mathring{g}} {\phi}+\lambda {\phi}=0\,\, &\text{in}\,\, C_{1+v}^T,\\
{\phi}=0 &\text{on}\,\, \partial C_{1+v}^T.
\end{array}
\right.
\end{equation}
It follows from \cite[Theorem 1.13]{Ambrosetti} that the problem (\ref{eigenvalueproblem}) possesses a sequence eigenvalues
\begin{equation}
0<\lambda_{1,v}<\lambda_{2,v}\leq\cdots, \lambda_{k,v}\nearrow+\infty.
\nonumber
\end{equation}

Let ${\phi}_{2,v}$ be an eigenfunction corresponding to $\lambda_{2,v}$ such that
\begin{equation}
\int_{C_{1+v}^{2 \pi}}{\phi}_{2,v}^2\left(x,\frac{T}{2 \pi}t\right)\,\text{dvol}_{\mathring{g}}=1,\nonumber
\end{equation}
where
\begin{equation}
\partial C_{1+v}^T=\left\{(x,t)\in \mathbb{R}^N\times \mathbb{R}/T\mathbb{Z}: \vert x\vert=1+v\left(\frac{2\pi t}{T}\right)\right\}.\nonumber
\end{equation}
Applying \cite[Theorem 11.4]{Gilbarg}, we see that $\phi_{2,v}\in \mathcal{C}^{2,\alpha}\left(\overline{C_{1+v}^T}\right)$.
Clearly, ${\phi}_{2,v}$ and $\lambda_{2,v}$ depend smoothly on $v$ (the smoothness can also be obtained by the Implicit Function Theorem as that of \cite[Proposition 4.1]{RSW}), and ${\phi}_{2,0}=\pm{\phi}_2$, $\lambda_{2,0}=\lambda_2$.
Without loss of generality, we assume ${\phi}_{2,0}={\phi}_2$.
For any fixed $t$, we have that ${\phi}_{2,v}$ is radially symmetric with respect to the first variable.
Hence, for any fixed $t$, ${\phi}_{2,v}(r,t)$ has a unique simple zero in $(0,1+v(t))$ which is denoted by $r_1(t)$.
The Implicit Function Theorem implies that $r_1(t)$ is $C^{2,\alpha}$ in local. By the arbitrariness of $t$, ${\phi}_{2,v}(x,t)$ is positive
in $C_{1+v}^{T,+}=\left\{(x,t)\in C_{1+v}^{T}:\vert x\vert<r_1(t)\right\}$ and is negative in $C_{1+v}^{T,-}=\left\{(x,t)\in C_{1+v}^{T}:\vert x\vert>r_1(t)\right\}$.
That is to say $r_1(t)$ is the unique zero line of ${\phi}_{2,v}(r,t)$.

Define the operator
\begin{equation}
\mathcal{N}(v,T)=\mathring{g}\left(\nabla {\phi}_v,\vartheta\right)\Big|_{\partial C_{1+v}^T}-\frac{1}{\text{Vol}_{\mathring{g}}\left(\partial C_{1+v}^T\right)}\int_{\partial C_{1+v}^T}\mathring{g}\left(\nabla {\phi}_v,\vartheta\right)\,\text{dvol}_{\mathring{g}},\nonumber
\end{equation}
where $\vartheta$ denotes the unit normal vector field to $\partial C_{1+v}^T$.
Note that $\mathcal{N}$ depends only on the variable $t$. Thus, we can define
\begin{equation}
F(v,T)(t)=\mathcal{N}(v,T)\left(\frac{T}{2\pi}t\right).\nonumber
\end{equation}
Since $\partial_r {\phi}_2(1)$ is a constant, it follows that $F(0,T)=0$ for any $T>0$.
Therefore, finding nontrivial domains emanating from $B_1\times \mathbb{R}$ such that the problem (\ref{overdeterminedeigenvalueproblem}) has a sign-changing solution
is equivalent to study the nontrivial solutions of $F(v,T)=0$.
In the spirit of Crandall-Rabinowitz local bifurcation theorem, it's key to find the degenerate point of linearization operator of $F$ and verify the transversality condition.
To study the linearization operator of $F$ with respect to $v$ at point $(0,T)$,
we first consider the following equation
\begin{equation}\label{ckequation1}
\left(\partial_r^2+\frac{N-1}{r}\partial_r+\lambda_2\right)c-\left(\frac{2k\pi}{T}\right)^2c=0
\end{equation}
with $c(1)=-\partial_r\phi_2(1)$ and $c'(0)=0$.
\\ \\
\textbf{Proposition 2.1.} \emph{For each $k\in \mathbb{N}$,
\begin{itemize}
  \item problem (\ref{ckequation1}) has a unique solution $c_k$ if $T\neq2k\pi/\sqrt{\lambda_2-\lambda_1}$,
  \item there is no solution to problem (\ref{ckequation1}) when $T=2k\pi/\sqrt{\lambda_2-\lambda_1}$.
\end{itemize}
Moreover, $c_k'$ is analytic if $T\neq2k\pi/\sqrt{\lambda_2-\lambda_1}$.}
\\ \\
\textbf{Proof.} In fact, for $T=2k\pi/\sqrt{\lambda_2-\lambda_1}$, we have that
\begin{equation}
\left(\partial_r^2+\frac{N-1}{r}\partial_r+\lambda_1\right)c=0.\nonumber
\end{equation}
Set $\vartheta(r):=c(r)+\partial_r\phi_2(1)$. Then we see that
\begin{equation}\label{vckequation}
\left(\partial_r^2+\frac{N-1}{r}\partial_r+\lambda_1\right)\vartheta=\lambda_1\phi_2'(1)
\end{equation}
with $\vartheta'(0)=0$ and $\vartheta(1)=0$.
Note that $\phi_1$ satisfies
\begin{equation}
\left\{
\begin{array}{ll}
\left(r^{N-1}\phi_1'\right)'+\lambda_1 r^{N-1}\phi_1=0\,\, &\text{in}\,\, (0,1),\\
\phi_1'(0)=\phi_1(1)=0.
\end{array}
\right.\nonumber
\end{equation}
Multiplying the equation by $\phi'_2(1)$ and integrating, we have that
\begin{equation}
\lambda_1\phi_2'(1)\int_0^1r^{N-1}\phi_1\,dr=-\phi_2'(1)\phi_1'(1)>0.\nonumber
\end{equation}
It follows that $\lambda_1\phi_2'(1)$ is not orthogonal to $\phi_1$.
By the Fredholm alternative theorem \cite[Theorem 6.2.4]{Evans}, problem (\ref{vckequation}) has no solution, which verifies the claim.

On the other hand, if $T\neq2k\pi/\sqrt{\lambda_2-\lambda_1}$,
we have that
\begin{equation}
\left(\partial_r^2+\frac{N-1}{r}\partial_r+\left(\lambda_2-\left(\frac{2k\pi}{T}\right)^2\right)\right)\vartheta=\left(\lambda_2-\left(\frac{2k\pi}{T}\right)^2\right)\phi_2'(1).\nonumber
\end{equation}
We see that
\begin{equation}
\lambda_2-\left(\frac{2k\pi}{T}\right)^2<\lambda_2\,\,\, \text{and}\,\,\,\lambda_2-\left(\frac{2k\pi}{T}\right)^2\neq\lambda_1 .\nonumber
\end{equation}
Using the Fredholm alternative theorem again, we see that problem (\ref{vckequation}) has a unique solution.
In conclusion, problem (\ref{ckequation1}) has a unique solution if and only if $T\neq2k\pi/\sqrt{\lambda_2-\lambda_1}$.
We use $c_k$ to denote this unique solution.

We finally prove the analyticity of $c_k'$. It is enough to show that $c_k$ is analytic with respect to $T$.
We use the following fact to show the analyticity of $c_k$:
if $F$ is an invertible operator, by the equality
\begin{equation}
(I-s F)^{-1}=\sum_{i\geq0}s^iF^i\nonumber
\end{equation}
for each $s\in \mathbb{R}$, the solution of
\begin{equation}
\left(F-\frac{\rho}{T^2} I\right)u=h\nonumber
\end{equation}
is analytic in $T$, where $I$ is the identity, $h$ is any continuous function and $\rho$ is a constant.

We first assume that $T>2k\pi/\sqrt{\lambda_2-\lambda_1}$. We consider
\begin{equation}
F=\partial_r^2+\frac{N-1}{r}\partial_r+\lambda_2 -\frac{1}{2}\left(\frac{2k\pi}{T}\right)^2\nonumber
\end{equation}
acting on $X:=\left\{c\in C[0,1]\cap C^2[0,1):c(1)=-\partial_r\phi_2(1), c'(0)=0\right\}$.
For any $f\in C[0,1]$, we consider
\begin{equation}
\left(\partial_r^2+\frac{N-1}{r}\partial_r+\lambda_2\right)c-\frac{1}{2}\left(\frac{2k\pi}{T}\right)^2c=f(x)\nonumber
\end{equation}
on $X$.
Then we see that
\begin{equation}\label{ckequationfx}
\left(\partial_r^2+\frac{N-1}{r}\partial_r+\left(\lambda_2-\frac{1}{2}\left(\frac{2k\pi}{T}\right)^2\right)\right)\vartheta=\left(\lambda_2-\frac{1}{2}\left(\frac{2k\pi}{T}\right)^2\right)\phi_2'(1)+f(x)
\end{equation}
with $\vartheta(1)=0$ and $\vartheta'(0)=0$.
Since $T>2k\pi/\sqrt{\lambda_2-\lambda_1}$, we have that
\begin{equation}
\lambda_1<\lambda_2-\left(\frac{2k\pi}{T}\right)^2<\lambda_2-\frac{1}{2}\left(\frac{2k\pi}{T}\right)^2<\lambda_2.\nonumber
\end{equation}
The Fredholm alternative theorem \cite[Theorem 6.2.5]{Evans} implies that there exists a unique solution of (\ref{ckequationfx}).
It follows that $F:X\longrightarrow C[0,1]$ is invertible.
Taking $h=0$ and $\rho=2k^2\pi^2$, the analyticity of $c_k$ is deduced.

We next consider the case of $T<2k\pi/\sqrt{\lambda_2-\lambda_1}$.
Now we take
\begin{equation}
F=\partial_r^2+\frac{N-1}{r}\partial_r+\lambda_2 -\left(\frac{2k\pi}{T}\right)^2.\nonumber
\end{equation}
For any $f\in C[0,1]$, similarly we have that
\begin{equation}\label{ckequationfx1}
\left(\partial_r^2+\frac{N-1}{r}\partial_r+\left(\lambda_2-\left(\frac{2k\pi}{T}\right)^2\right)\right)\vartheta=\left(\lambda_2-\left(\frac{2k\pi}{T}\right)^2\right)\phi_2'(1)+f(x)
\end{equation}
with $\vartheta(1)=0$ and $\vartheta'(0)=0$. The fact $T<2k\pi/\sqrt{\lambda_2-\lambda_1}$ implies that
\begin{equation}
\lambda_2-\left(\frac{2k\pi}{T}\right)^2<\lambda_1.\nonumber
\end{equation}
Using the Fredholm alternative theorem \cite[Theorem 6.2.5]{Evans} again, we deduce the existence and uniqueness of the solution of (\ref{ckequationfx1}).
Hence $F:X\longrightarrow C[0,1]$ is still invertible. Taking $h=0$ and $\rho=0$, the analyticity of $c_k$ is concluded again.
\qed\\

By Fourier expansion $v$ can be written as
\begin{equation}
v=\sum_{k\geq1}a_k\cos(kt).\nonumber
\end{equation}
Denote
\begin{equation}
\sum_{k\geq1}c_k(r)a_k\cos\left(\frac{2k\pi t}{T}\right):=\psi.\nonumber
\end{equation}
Then we can verify that $\psi$ is $L^2\left(C_1^T\right)$-orthogonal to $\phi_2$ and
satisfies the following problem
\begin{equation}\label{eigenvalueonc1cvnu=0v=0}
\left\{
\begin{array}{ll}
\Delta_{\mathring{g}}  \varphi+\lambda_2\varphi=0\,\, &\text{in}\,\, C_{1}^{T},\\
\varphi=-\partial_r{\phi}_2v\left(\frac{2\pi t}{T}\right) &\text{on}\,\, \partial C_{1}^{T}.
\end{array}
\right.
\end{equation}
Conversely, by Proposition 2.1, we deduce that, for $T\neq2k\pi/\sqrt{\lambda_2-\lambda_1}$, the problem (\ref{eigenvalueonc1cvnu=0v=0}) has a unique solution which is just $\psi$.

Form now on we always assume that $T\neq2k\pi/\sqrt{\lambda_2-\lambda_1}$ for each $k\in \mathbb{N}$.
Define
\begin{equation}
\widetilde{H}_T(v)=\left(\partial_r \psi+\partial_r^2\phi_2 v\left(\frac{2\pi}{T}t\right)\right)\bigg|_{\partial C_{1}^{T}}\nonumber
\end{equation}
and
\begin{equation}
{H}_T(v)=\widetilde{H}_T(v)\left(\frac{T}{2\pi}t\right).\nonumber
\end{equation}
By an argument as that of \cite[Proposition 3.4]{Sicbaldi} with obvious changes, the linearized operator of $F$ with respect to $v$ at point $(0,T)$ is just
${H}_T$.

Let $V_k$ be the space spanned by the function $\cos(kt)$.
The variable separation characteristics of $\psi$ shows that $H_T$ preserves the eigenspace $V_k$.
Let $\sigma_k(T)$ be the eigenvalues of $H_T$ with respect to the eigenfunctions $\cos(kt)$.
Similar to that \cite{Sicbaldi} we have that
\begin{equation}
\sigma_k(T)=\partial_r c_k(1)+\partial_r^2 \phi_2(1),\nonumber
\end{equation}
where $c_k$ is the continuous solution on $[0,1]$ of (\ref{ckequation1}).
Note that $\sigma_k(T)=\sigma_1\left(T/k\right)$, which indicates that the property of $\sigma_k$ can be deduced from the property of $\sigma_1$. We next only consider the case of $k=1$.
The analyticity of $c_1'$ implies that $\sigma_1(T)$ is differentiable.
The zero (if it exists) of $\sigma_1(T)$ is just the degenerate point of ${H}_T$.
We postpone the study of the zero of $\sigma_1(T)$ in the next section.
We end this section by presenting some conclusions on Bessel functions (we refer to \cite[Chapter 10]{Olver} for details), which will be used later.

\subsection{Some conclusions on Bessel functions}

\quad\, For $\tau\in \mathbb{R}$, the Bessel function of the first kind is defined by
\begin{equation}
J_\tau(s)=\sum_{m=0}^\infty\frac{(-1)^m\left(\frac{s}{2}\right)^{2m+\tau}}{m!\Gamma(\tau+m+1)},\nonumber
\end{equation}
which is the solution of the differential equation
\begin{equation}
s^2f''(s)+sf'(s)+\left(s^2-\tau^2\right)f(s)=0.\nonumber
\end{equation}
If $\tau$ is an integer, then $J_{-\tau}(s)=(-1)^\tau J_{\tau}(s)$.
For any $\tau\in \mathbb{R}$ and $s>0$ we have the following relations
\begin{equation}\label{relationsforbesselfunction1}
J_{\tau-1}(s)-J_{\tau+1}(s)=2J_{\tau}'(s),
\end{equation}
\begin{equation}\label{relationsforbesselfunction2}
sJ_{\tau}'(s)+\tau J_{\tau}(s)=sJ_{\tau-1}(s)
\end{equation}
and
\begin{equation}\label{relationsforbesselfunction3}
sJ_{\tau}'(s)-\tau J_{\tau}(s)=-sJ_{\tau+1}(s).
\end{equation}
When $\tau$ is fixed and $s\rightarrow 0$, it has that
\begin{equation}\label{relationsforbesselfunctionz=0}
J_\tau(s)\thicksim\frac{\left(\frac{s}{2}\right)^\tau}{\Gamma(\tau+1)}
\end{equation}
for $\tau\neq-1,-2,-3,\ldots$. In addition, one has that $J_0(0)=1$ and $J_\tau(0)=0$ for all $\tau>0$.

It is well known that the eigenvalues and eigenfunctions of the Dirichlet Laplacian on the unit ball have the following relations
\begin{equation}
\lambda_k=j_{\nu,k}^2\,\,\,\text{and}\,\,\,\vert x\vert^{\nu}\overline{\phi}_k=C_kJ_\nu\left(j_{\nu,k}\vert x\vert\right),\nonumber
\end{equation}
where $j_{\nu,k}$ is the $k$-th positive zero of $J_{\nu}$ for $\nu=(N-2)/2$ with $N\geq2$. Denote
$$
\mu=:\frac{2\pi}{j_{\nu,2}}=\frac{2\pi}{\sqrt{\lambda_2}}.
$$
Here we assume that $\overline{\phi}_1$ is positive, then $J_\nu$ is positive in $\left(0,\sqrt{\lambda_1}\right)$.

We claim that the above relation is also valid for $N=1$.
We know that $\overline{\phi}_k$ satisfies
\begin{equation}
\phi''+\lambda_k\phi=0\nonumber
\end{equation}
with $\phi(1)=0$ and $\phi'(0)=0$.
Then we can verify that $\widehat{\phi}_k(s)=s^\tau \overline{\phi}_k\left(s/\sqrt{\lambda_k}\right)$ satisfies
\begin{equation}
s^2\widehat{\phi}_k''(s)+s\widehat{\phi}_k'(s)+\left(s^2-\tau^2\right)\widehat{\phi}_k(s)=0,\nonumber
\end{equation}
where $s=\sqrt{\lambda_k}r$ and $\tau=-1/2$.
It follows that
\begin{equation}
\widehat{\phi}_k(s)=C_kJ_{-\frac{1}{2}}(s).\nonumber
\end{equation}
Further,
\begin{equation}
\phi\left(r\right)=\left(\sqrt{\lambda_k}r\right)^{\frac{1}{2}}C_kJ_{-\frac{1}{2}}(\sqrt{\lambda_k}r).\nonumber
\end{equation}
In particular,
\begin{equation}\label{J12}
J_{-\frac{1}{2}}\left(\frac{(2k-1)\pi}{2}r\right)=\left(\frac{(2k-1)\pi}{2}r\right)^{-\frac{1}{2}}\cos\left(\frac{(2k-1)\pi}{2}r\right),
\end{equation}
where we take $C_k=1$ for simplicity. Note that formula (\ref{J12}) also can be derived by the following formula (see \cite[Chapter 10, formula 10.16.1]{Olver})
\begin{equation}
J_{-\frac{1}{2}}(z)=\left(\frac{2}{\pi z}\right)^{\frac{1}{2}}\cos z.\nonumber
\end{equation}
This implies that $\lim_{z\rightarrow0^+}J_{-\frac{1}{2}}(z)=+\infty$.

Since there exists an explicit expression of $\overline{\phi}_k$ when $N=3$, there should be an explicit expression of $J_{\frac{1}{2}}(z)$.
Now $\overline{\phi}_k$ satisfies
\begin{equation}
\phi''+\frac{2}{r}\phi'+\lambda_k\phi=0\nonumber
\end{equation}
with $\phi(1)=0$ and $\phi'(0)=0$.
Since $\tau=1/2$, we can verify that $\widehat{\phi}_k(s)=s^\tau \overline{\phi}_k\left(s/\sqrt{\lambda_k}\right)$ and it satisfies
\begin{equation}
s^2\widehat{\phi}_k''(s)+s\widehat{\phi}_k'(s)+\left(s^2-\tau^2\right)\widehat{\phi}_k(s)=0,\nonumber
\end{equation}
which implies that
\begin{equation}
\widehat{\phi}_k(s)=C_kJ_{\frac{1}{2}}(s).\nonumber
\end{equation}
Further,
\begin{equation}
\phi\left(r\right)=\left(\sqrt{\lambda_k}r\right)^{-\frac{1}{2}}C_kJ_{\frac{1}{2}}(\sqrt{\lambda_k}r).\nonumber
\end{equation}
In particular, for $r>0$, using Lemma 2.1, we obtain that
\begin{equation}\label{J123}
J_{\frac{1}{2}}\left(k\pi r\right)=\left(k\pi r\right)^{\frac{1}{2}}\frac{\sin(k\pi r)}{r},
\end{equation}
up to a positive constant factor. This is consistent with the following formula (see \cite[Chapter 10, formala 10.16.1]{Olver})
\begin{equation}
J_{\frac{1}{2}}(z)=\left(\frac{2}{\pi z}\right)^{\frac{1}{2}}\sin z,\nonumber
\end{equation}
up to a positive constant factor. In particular, the first positive zero of $J_{\frac{1}{2}}$ is $\pi$ and the second positive zero of $J_{\frac{1}{2}}$ is $2\pi$.

If $\tau$ is real, then $J_\tau(s)$ has an infinite number of positive real zeros. All of
these zeros are simple, the $m$-th positive zero of $J_\tau(s)$ is denoted by $j_{\tau,m}$.
When $\tau\geq0$, the zeros interlace according to the
inequalities
\begin{equation}\label{interlaceproperty}
j_{\tau,1}< j_{\tau+1,1} < j_{\tau,2} < j_{\tau+1,2} < j_{\tau,3} <\cdots.
\end{equation}
When $\tau\geq-1$ the zeros of $J_\tau(s)$ are all real.

For $\tau\in \mathbb{R}$, the modified Bessel function of the first kind is defined by
\begin{equation}\label{secondbesselfunction}
I_\tau(s)=\sum_{m=0}^\infty\frac{\left(\frac{s}{2}\right)^{2m+\tau}}{m!\Gamma(\tau+m+1)},
\end{equation}
which is the solution of
\begin{equation}
s^2f''(s)+sf'(s)-\left(s^2+\tau^2\right)f(s)=0.\nonumber
\end{equation}

\section{The properties of $\sigma_1(T)$ for $N\geq 2$}

\quad\, In this section, we prove the existence of zero to $\sigma_1(T)$ and study its properties at zero.
In particular, we would establish that the derivatives of $\sigma_1$ at its zero (if it exists) does not vanish because it is the key to verify the transversality condition.
For $N\geq2$, we first obtain the asymptotic behavior of $\sigma_1$ as follows.
\\ \\
\textbf{Proposition 3.1.} \emph{Let
\begin{equation}
\mu=\frac{2\pi}{\sqrt{\lambda_2}}<\frac{2\pi}{\sqrt{\lambda_2-\lambda_1}}=\delta.\nonumber
\end{equation}
The function $\sigma(T):=\sigma_1(T)$ is well defined on $\mathbb{R}\setminus{\delta}$, and $T=\delta$ is its only singular point.
Moreover, $\sigma(T)$ satisfies
$\lim_{T\rightarrow 0^+}\sigma(T)=-\infty$, $\lim_{T\rightarrow \delta^-}\sigma(T)=+\infty$, $\lim_{T\rightarrow \delta^+}\sigma(T)=-\infty$ and $\lim_{T\rightarrow+\infty}\sigma(T)=+\infty$}.
\\ \\
\textbf{Proof.} 
It is well known that
$\sqrt{\lambda_2}$ is the second zero of the first kind Bessel function $J_\nu$.
For $T\in(0,\mu)$, by the computations of \cite[Section 5]{Schlenk} with obvious changes we have that
\begin{equation}
\sigma(T)=-\kappa_n\sqrt{\lambda_2}^{1-\nu}J'_\nu\left(\sqrt{\lambda_2}\right)\left(1+\frac{\xi I_{\nu-1}(\xi)}{I_{\nu}(\xi)}\right),\nonumber
\end{equation}
where $\kappa_n$ is a positive constant, $I_{\nu}$ is the first kind modified Bessel function, $\xi=\sqrt{\left(\frac{2\pi}{T}\right)^2-\lambda_2}$.
Based on the fact $\lim_{s\rightarrow\infty}\frac{I_{\nu}(s)}{\frac{1}{\sqrt{2\pi s}}e^s}=1$,
we have that
\begin{equation}
\lim_{T\rightarrow0}\frac{\xi I_{\nu-1}(\xi)}{I_{\nu}(\xi)}=\lim_{T\rightarrow0}\xi=\infty.\nonumber
\end{equation}
Thus, it follows that
\begin{equation}
\lim_{T\rightarrow0}\sigma(T)=-\infty, \nonumber
\end{equation}
where we use $J'_\nu\left(\sqrt{\lambda_2}\right)>0$.

For $T>\mu$, we similarly get
\begin{eqnarray}\label{formaright}
\sigma(T)&=&-\kappa_n\sqrt{\lambda_2}^{1-\nu}J'_\nu\left(\sqrt{\lambda_2}\right)\left(1+\frac{\rho J_{\nu-1}(\rho)}{J_{\nu}(\rho)}\right)\nonumber\\
&=&-\kappa_n\sqrt{\lambda_2}^{1-\nu}J'_\nu\left(\sqrt{\lambda_2}\right)\left((2\nu+1)-\frac{\rho J_{\nu+1}(\rho)}{J_{\nu}(\rho)}\right),
\end{eqnarray}
where
\begin{equation}
\rho=\sqrt{\lambda_2-\left(\frac{2\pi}{T}\right)^2}.\nonumber
\end{equation}
We clearly see $\rho\nearrow\sqrt{\lambda_2}$ as $T\rightarrow+\infty$.
When $T=\delta$, we see that $\rho=\sqrt{\lambda_1}$ is the only zero of $J_{\nu}$ in $\left(0,\sqrt{\lambda_2}\right)$.
From (\ref{formaright}) we get that $T=\delta$ is the only singular point of $\sigma$.
Since $\rho J_{\nu+1}(\rho)$ and $J_{\nu}(\rho)$ have the same (opposite) sign on the left (right) side of $\rho=\sqrt{\lambda_1}$ and $J'_\nu\left(\sqrt{\lambda_2}\right)>0$, we obtain that
$$\lim_{T\rightarrow \delta^-}\sigma(T)=+\infty,\quad \lim_{T\rightarrow \delta^+}\sigma(T)=-\infty,$$
where we use (\ref{formaright}).

Since $j_{\nu,2}<j_{\nu+1,2}$, we see that $J_{\nu+1}\left(\sqrt{\lambda_2}\right)<0$. Then, as the argument of \cite[Lemma 6.1]{Schlenk}, we have that
\begin{equation}
\lim_{\rho\nearrow \sqrt{\lambda_2}}\frac{\rho J_{\nu+1}(\rho)}{J_{\nu}(\rho)}=+\infty.\nonumber
\end{equation}
Since $J'_\nu\left(\sqrt{\lambda_2}\right)>0$, we obtain
\begin{equation}
\lim_{T\rightarrow\infty}\sigma(T)=\lim_{\rho\nearrow \sqrt{\lambda_2}}\frac{\rho J_{\nu+1}(\rho)}{J_{\nu}(\rho)}=+\infty\nonumber
\end{equation}
as desired.\qed\\

We next study the differentiability and monotonicity of $\sigma(T)$.
\\ \\
\textbf{Proposition 3.2.} \emph{The function $\sigma(T)$ is strictly increasing on $(0,\delta)$ and $(\delta,+\infty)$, respectively.
Moreover, $\sigma(T)<0$ for $T\leq\mu$.}
\\ \\
\textbf{Proof.} For $T<\mu$, as that of \cite[Section 5]{Schlenk} we can show that
\begin{equation}
\sigma(T)=\varphi_2''(1)-\varphi_2'(1)\frac{\xi I_{\nu+1}(\xi)}{I_{\nu}(\xi)}.\nonumber
\end{equation}
Let
\begin{equation}
f(s)=\frac{s I_{\nu+1}(s)}{I_{\nu}(s)}.\nonumber
\end{equation}
From the argument of \cite[Lemma 6.3]{Schlenk}, we know that $f'(s)>0$ for all $s\in(0,+\infty)$.
Since $\xi'(T)<0$ for all $T\in(0,\mu)$, we have that $f(\xi(T))$ is differentiable with respect to $T$ and
\begin{equation}
\frac{\text{d}}{\text{d} T}f(\xi(T))<0.\nonumber
\end{equation}
Since $\varphi_2'(1)>0$, we obtain that $\sigma'(T)>0$.

By (\ref{secondbesselfunction}) we have that
\begin{equation}
\lim_{\xi\rightarrow0}\frac{\xi I_{\nu+1}(\xi)}{I_{\nu}(\xi)}=0\,\,\,\text{or}\,\,\,\lim_{\xi\rightarrow0}\frac{\xi I_{\nu-1}(\xi)}{I_{\nu}(\xi)}=2\frac{\Gamma(\nu+1)}{\Gamma(\nu)}=2\nu.\nonumber
\end{equation}
It is worth noting that there may be a typo when computing the first limit in \cite[Lemma 5.1]{Schlenk}.
Then, as that of \cite[Lemma 5.1]{Schlenk} we can show that
\begin{equation}
\sigma(\mu)=-\kappa_n\sqrt{\lambda_2}^{1-\nu}J'_\nu\left(\sqrt{\lambda_2}\right)\left(2\nu+1\right)<0\nonumber
\end{equation}
due to $J'_\nu\left(\sqrt{\lambda_2}\right)>0$ and $2\nu+1=N-1>0$.
Therefore, $\sigma(T)<0$ for $T\leq\mu$.


We now consider the case of $T>\mu$. We divide the range of $T$ into the following two subintervals
\begin{equation}
T<\delta\nonumber
\end{equation}
and
\begin{equation}
T>\delta.\nonumber
\end{equation}
We first consider the range of $T\in\left(\mu,\delta\right)$.
From the definition of $\rho(T)$ we have that
\begin{equation}
\rho:\left(\mu,\delta\right)\rightarrow\left(0,\sqrt{\lambda_1}\right)\nonumber
\end{equation}
and it is strictly increasing.
Set
\begin{equation}
h(s)=\frac{s J_{\nu+1}(s)}{J_{\nu}(s)}.\nonumber
\end{equation}
Then we see that the sign of $h'(s)$ is determined by the sign of $J_\nu^2(s)-J_{\nu-1}(s)J_{\nu+1}(s)$.
From the argument of \cite[Claim 6.6]{Schlenk} (in fact, there is a gap, which will be corrected in Appendix) we know that
$\sigma'(T)$ in $\left(\mu,\delta\right)$ has the same sign as $h'(s)$ in $s\in \left(0,\sqrt{\lambda_1}\right)$.
From Claim A.1 in Appendix, we know that, $J_\nu^2(s)>J_{\nu-1}(s)J_{\nu+1}(s)$ for all $s\in\left(0,\sqrt{\lambda_1}\right)$.
This combines (\ref{formaright}) and
 $J'_\nu\left(\sqrt{\lambda_2}\right)>0$ implies that
\begin{equation}
\sigma'(T)>0\nonumber
\end{equation}
for $T\in\left(\mu,\delta\right)$.

We claim that $\sigma'(\mu)>0$.
For $T<\mu$, since $\xi'(T)=-4\pi^2/T^3\xi(T)$, we have that
\begin{eqnarray}
\sigma'(T)&=&-\varphi_2'(1)f'(\xi(T))\xi'(T)\nonumber\\
&=&-\varphi_2'(1)\frac{\xi(T)\left(I_\nu^2(\xi(T))-I_{\nu-1}(\xi(T))I_{\nu+1}(\xi(T))\right)}{I_\nu^2(\xi(T))}\xi'(T)\nonumber\\
&=&4\pi^2\varphi_2'(1)\frac{\left(I_\nu^2(\xi(T))-I_{\nu-1}(\xi(T))I_{\nu+1}(\xi(T))\right)}{T^3I_\nu^2(\xi(T))}.\nonumber
\end{eqnarray}
Since $\lim_{T\nearrow\mu}\xi(T)=0$, the left derivative of $\sigma$ at $T=\mu$ is equal to
\begin{equation}
\lim_{T\nearrow\mu}\sigma'(T)=\frac{4\pi^2\varphi_2'(1)}{\mu^3}\left(1-\lim_{s\searrow0}\frac{I_{\nu-1}(s)I_{\nu+1}(s)}{I_\nu^2(s)}\right).\nonumber
\end{equation}
It follows the power series expansion of $I_\nu$ that
\begin{equation}
I_\nu(s)=\frac{\left(\frac{s}{2}\right)^{\nu}}{\Gamma(\nu+1)}+O\left(s^{\nu+2}\right).\nonumber
\end{equation}
Then we compute that
\begin{equation}
\lim_{s\searrow0}\frac{I_{\nu-1}(s)I_{\nu+1}(s)}{I_\nu^2(s)}=\frac{\Gamma^2(\nu+1)}{\Gamma(\nu)\Gamma(\nu+2)}=\frac{\nu}{\nu+1}<1.\nonumber
\end{equation}
Therefore, we obtain that
the left derivative of $\sigma$ at $T=\mu$ is positive and equal to
\begin{equation}
\lim_{T\nearrow\mu}\sigma'(T)=\frac{4\pi^2\varphi_2'(1)}{\mu^3}\left(1-\frac{\Gamma^2(\nu+1)}{\Gamma(\nu)\Gamma(\nu+2)}\right).\nonumber
\end{equation}
On the other hand, when $T>\mu$, as that of \cite[Section 5]{Schlenk} we can show that
\begin{equation}\label{rightsigma}
\sigma(T)=\varphi_2''(1)+\varphi_2'(1)\frac{\rho J_{\nu+1}(\rho)}{J_{\nu}(\rho)}.
\end{equation}
Then, reasoning as the left derivative, we obtain that the right derivative $\sigma$ at $T=\mu$ is positive and equal to
\begin{equation}
\lim_{T\searrow\mu}\sigma'(T)=\frac{4\pi^2\varphi_2'(1)}{\mu^3}\left(1-\frac{\Gamma^2(\nu+1)}{\Gamma(\nu)\Gamma(\nu+2)}\right)=\lim_{T\nearrow\mu}\sigma'(T),\nonumber
\end{equation}
which is the desired claim.
Of course, the differentiability at $T=\mu$ is also a consequence of analyticity.

We further consider the case of $T>\delta$.
For this case we see that $\rho(T)\in\left(\sqrt{\lambda_1}, \sqrt{\lambda_2}\right)$ and $\rho$ is differentiable and strictly increasing on $[\delta,+\infty)$.
Note that the interlace property (\ref{interlaceproperty}) may not be valid for $\tau<0$ (see the Appendix).
While, $\nu-1$ is negative when $N=2$ or $N=3$. The interlace property (\ref{interlaceproperty}) cannot be used directly to deal with the case of $N=2,3$.

We first consider the case of $N\geq4$. In this case, the zeros interlace for $\nu-1\geq0$ and the profiles of $J(s)$ are as follows,
then we have that
\begin{equation}
j_{\nu-1,2},\,\, j_{\nu+1,1}\in\left(\sqrt{\lambda_1},\sqrt{\lambda_2}\right).\nonumber
\end{equation}
Set
\begin{equation}
\alpha=\min\left\{j_{\nu-1,2},j_{\nu+1,1}\right\},\,\, \beta=\max\left\{j_{\nu-1,2},j_{\nu+1,1}\right\}.\nonumber
\end{equation}
One of the following two cases will occur:
$$
(a)~\alpha<\beta; \quad (b)~\alpha=\beta.
$$
\begin{figure}[ht]
\centering
\includegraphics[width=1\textwidth]{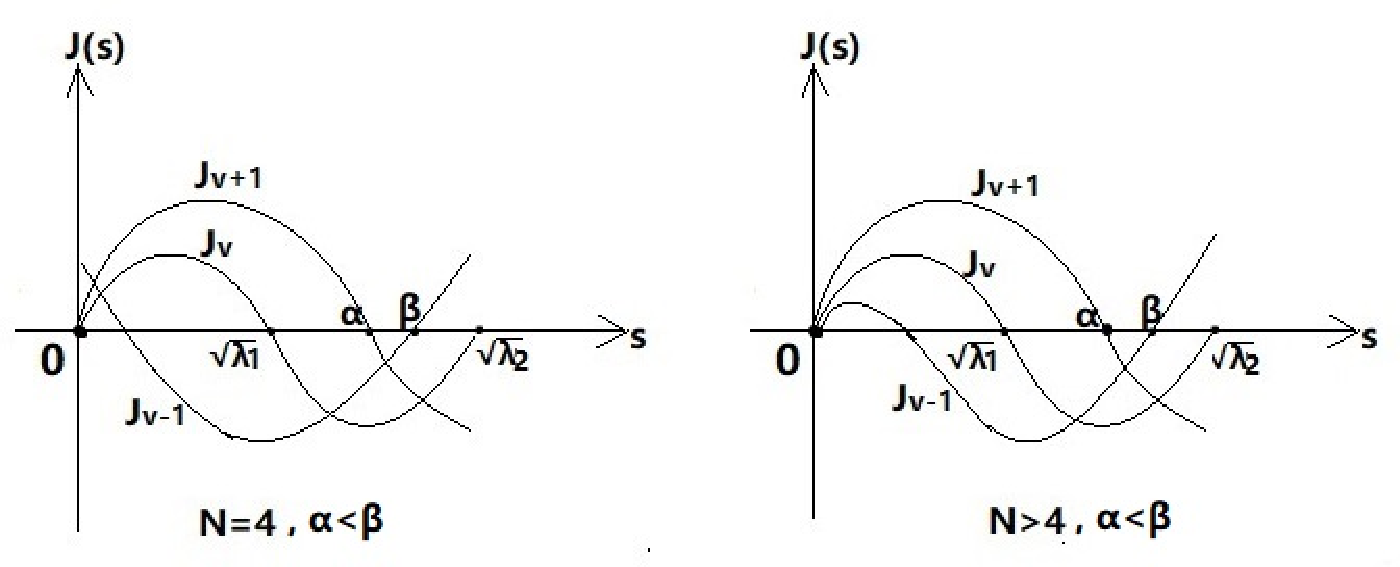}
\caption{The profile of $J(s)$ for $N\geq 4$. }
\end{figure}
\indent We next only prove the case of (a) because the case of (b) is similar.
For $s\in\left(\sqrt{\lambda_1},\alpha\right]$, we know that
\begin{equation}
J_{\nu-1}(s)\leq0, \,\,J_{\nu+1}(s)\geq0.\nonumber
\end{equation}
It follows that
\begin{equation}
J_{\nu-1}(s)J_{\nu+1}(s)\leq0.\nonumber
\end{equation}
Note that
\begin{equation}
h'(s)=\frac{s\left(J_\nu^2-J_{\nu-1}J_{\nu+1}\right)}{J_\nu^2}\,\,\,\text{for}\,\,\,s\in\left(\sqrt{\lambda_1}, \sqrt{\lambda_2}\right).\nonumber
\end{equation}
Thus, $h'(s)>0$ for $s\in\left(\sqrt{\lambda_1},\alpha\right]$.
Similarly, we can show that $h'(s)>0$ for $s\in\left[\beta,\sqrt{\lambda_2}\right)$.

We now assume $s\in\left(\alpha,\beta\right)$. 
In this interval, it is easy to check that $J_{\nu-1}(s)$ and $J_{\nu+1}(s)$ have the same signs (see Figure 1).

Thus we have that
\begin{equation}\label{sighaiphabelta}
-2J_{\nu-1}(s)J_{\nu+1}(s)<0
\end{equation}
for $s\in\left(\alpha,\beta\right)$.
From the relations (\ref{relationsforbesselfunction1})--(\ref{relationsforbesselfunction3})
we deduce that
\begin{equation}
2sJ_{\nu}'(s)=s\left(J_{\nu-1}(s)-J_{\nu+1}(s)\right),\nonumber
\end{equation}
\begin{equation}
sJ_{\nu-1}'(s)=(\nu-1)J_{\nu-1}(s)-sJ_{\nu}(s)\nonumber
\end{equation}
and
\begin{equation}
sJ_{\nu+1}'(s)=sJ_{\nu}(s)-(\nu+1)J_{\nu+1}(s).\nonumber
\end{equation}
Using these relations, we can verify that inequality (\ref{sighaiphabelta}) is equivalent to
\begin{equation}
2J_{\nu}(s)sJ_{\nu}'(s)>sJ_{\nu-1}'(s)J_{\nu+1}(s)+J_{\nu-1}(s)sJ_{\nu+1}'(s).\nonumber
\end{equation}
That is to say that
\begin{equation}
\left(J_{\nu}^2-J_{\nu-1}J_{\nu+1}\right)'(s)>0.\nonumber
\end{equation}
It follows that $J_{\nu}^2(s)-J_{\nu-1}(s)J_{\nu+1}(s)$ is strictly increasing in $\left(\alpha,\beta\right)$.
Since $J_{\nu}(\alpha)<0$ and $J_{\nu-1}(\alpha)J_{\nu+1}(\alpha)=0$, we have $J_{\nu}^2(\alpha)-J_{\nu-1}(\alpha)J_{\nu+1}(\alpha)>0$.
Thus we get that
\begin{equation}
J_{\nu}^2(s)-J_{\nu-1}(s)J_{\nu+1}(s)>0\nonumber
\end{equation}
for $s\in\left[\alpha,\beta\right]$. It follows that $h'(s)>0$ for $s\in\left[\alpha,\beta\right]$.
In conclusion, when $N\geq4$, we obtain that
$h'(s)>0$ for $s\in\left(\sqrt{\lambda_1},\sqrt{\lambda_2}\right)$.
In view of (\ref{rightsigma}) we conclude that $\sigma'(T)>0$ for $T>\delta$ and $N\geq4$.

We now consider the case of $N=2$.
In this case, $\nu=0$, $J_{-1}(s)=-J_1(s)$ and the profiles of $J(s)$ are as follows.
It follows that $J_{\nu-1}(s)J_{\nu+1}(s)=-J_1^2(s)\leq0$.
Hence we have that $J_{\nu}^2(s)-J_{\nu-1}(s)J_{\nu+1}(s)>0$ in $\left(\sqrt{\lambda_1},\sqrt{\lambda_2}\right)$.
We still have that $h'(s)>0$ for $s\in\left(\sqrt{\lambda_1},\sqrt{\lambda_2}\right)$.
\begin{figure}[ht]
\centering
\includegraphics[width=1\textwidth]{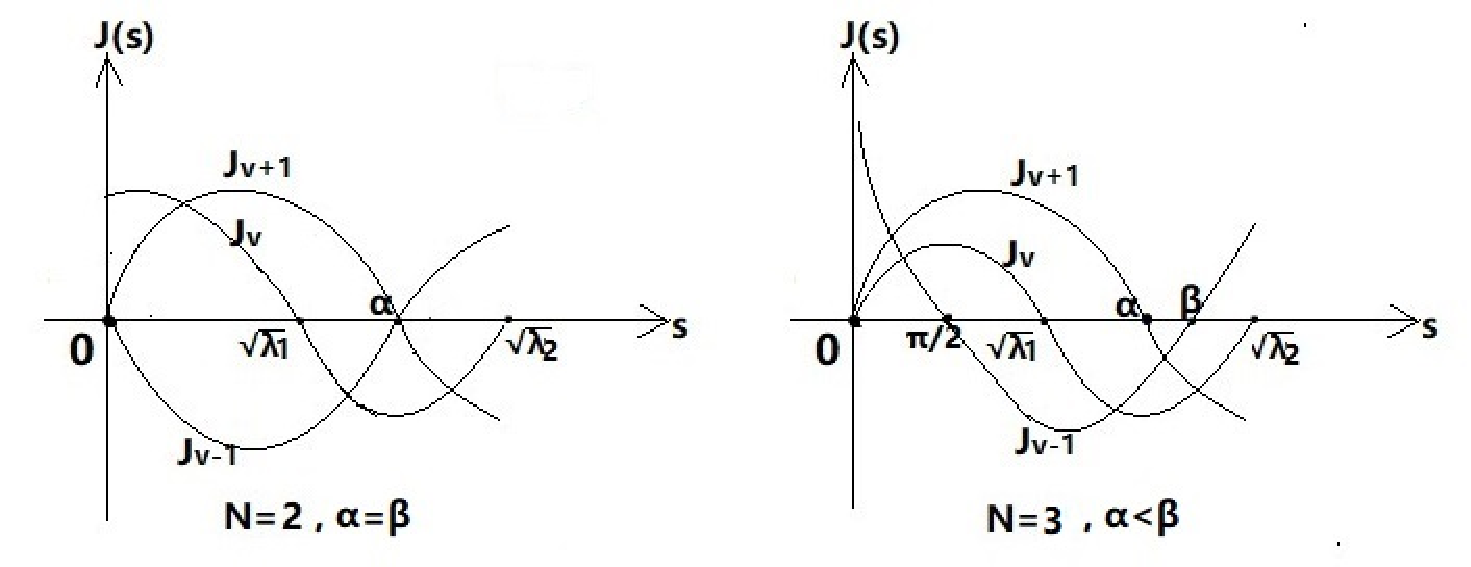}
\caption{The profile of $J(s)$ for $N=2$ and $N=3$. }
\end{figure}

We finally study the case of $N=3$. In this case, $\nu=1/2$, $\sqrt{\lambda_1}=\pi$ and $\sqrt{\lambda_2}=2\pi$ due to Lemma 2.1 and (\ref{J123}).
From (\ref{J12}) we see that
$J_{-1/2}(s)$ is positive in $\left(0,\pi/2\right)$, the first and second positive zeros are $\pi/2$ and $3\pi/2$.
We now define
\begin{equation}
\alpha=\min\left\{\frac{3\pi}{2},j_{\frac{3}{2},1}\right\},\,\, \beta=\max\left\{\frac{3\pi}{2},j_{\frac{3}{2},1}\right\}.\nonumber
\end{equation}
Without loss of generality, we still assume that $\alpha<\beta$.
We note that $J_{-1/2}(s)$ and $J_{3/2}(s)$ have
opposite signs in $\left(\sqrt{\lambda_1},\alpha\right]\cup \left[\beta,\sqrt{\lambda_2}\right)$ (see Figure 2).
So, we have that $h'(s)>0$ for $s\in\left(\sqrt{\lambda_1},\alpha\right]\cup \left[\beta,\sqrt{\lambda_2}\right)$.
Reasoning as in the case of $N\geq4$ we see that $J_{\frac{1}{2}}^2(s)-J_{-\frac{1}{2}}(s)J_{\frac{3}{2}}(s)$ is strictly increasing in $\left(\alpha,\beta\right)$.
Since $J_{-\frac{1}{2}}(\alpha)J_{\frac{3}{2}}(\alpha)=0$ and $J_{\frac{1}{2}}(\alpha)<0$, we have that $J_{\frac{1}{2}}^2(\alpha)-J_{-\frac{1}{2}}(\alpha)J_{\frac{3}{2}}(\alpha)>0$.
Thus we obtain that
\begin{equation}
J_{\frac{1}{2}}^2(s)-J_{-\frac{1}{2}}(s)J_{\frac{3}{2}}(s)>0\nonumber
\end{equation}
for $s\in\left[\alpha,\beta\right]$. It follows that $h'(s)>0$ for $s\in\left[\alpha,\beta\right]$. Therefore, we obtain that
$h'(s)>0$ for $s\in\left(\sqrt{\lambda_1},\sqrt{\lambda_2}\right)$ for $N=3$.
\qed\\

From Proposition 3.2 we conclude that $\sigma(T)$ has exactly two zeros $T_*$ and $T^*$ such that
\begin{equation}
T_*\in(\mu,\delta)\nonumber
\end{equation}
and
\begin{equation}
T^*>\delta.\nonumber
\end{equation}
Moreover, one has that $\sigma'\left(T_*\right)>0$ and $\sigma'\left(T^*\right)>0$.
Summing up, we clearly obtain the behaviors of eigenvalue $\sigma(T)$ as follows.
\begin{figure}[ht]
\centering
\includegraphics[width=0.8\textwidth]{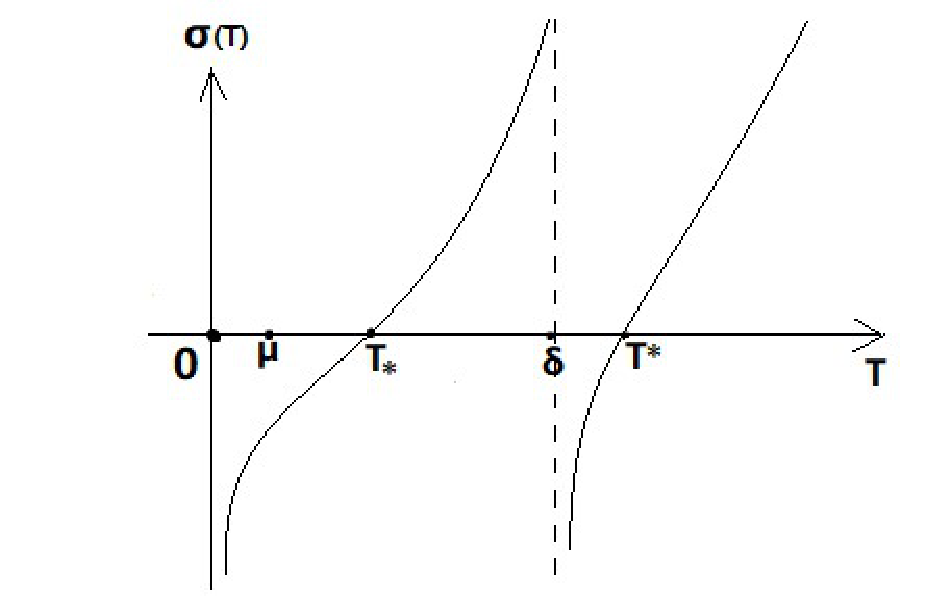}
\caption{The profile of $\sigma(T)$. }
\end{figure}

Thus, we further have the following result.
\\
\\
\textbf{Proposition 3.3.} \emph{The kernel space of $D_v F\left(0,T_*\right)$ is just $V_1$, where $V_1$ is the space spanned by the function $\cos(t)$. Moreover, if $T^*\neq iT_*$ for any $i\in \mathbb{N}$, the kernel space of $D_v F\left(0,T^*\right)$ is still $V_1$. While, if there exists some $i\in \mathbb{N}$ such that $T^*=i T_*$, the kernel space of $D_v F\left(0,T^*\right)$ is just $V_1\cup V_i$.}
\\ \\
\textbf{Proof.} Clearly, $V_1$ is contained in the kernel space of $D_v F\left(0,T_*\right)$ or $D_v F\left(0,T^*\right)$.
Note that
\begin{equation}
\sigma_{k}(T)=\sigma\left(\frac{T}{k}\right)\nonumber
\end{equation}
for all $k\in\mathbb{N}$.
As $\sigma\left(T_*/j\right)<0$ with any $j\geq2$, it follows that $V_j$ is not the kernel $D_v F\left(0,T_*\right)$. Thus, the kernel space of $D_v F\left(0,T_*\right)$ is just $V_1$.

If $T^*\neq iT_*$ for any $i\in \mathbb{N}$, then $\sigma\left(T^*/i\right)\neq0$ for any $i\geq2$. Hence $V_i$ is not the kernel $D_v F\left(0,T^*\right)$ for any $i\geq2$.
So, the kernel space of $D_v F\left(0,T^*\right)$ is still $V_1$. If there exists some $i\in \mathbb{N}$ with $i\geq2$ such that $T^*=i T_*$,
we have that $\sigma\left(T^*/i\right)=\sigma\left(T_*\right)=0$. Hence $V_i$ is also contained in the kernel space of $D_v F\left(0,T^*\right)$.
In this case, of course, $T^*\neq j T_*$ for any $j\neq i$ with $j\geq2$. Thus, $V_j$ is not contained the kernel $D_v F\left(0,T^*\right)$ for any $j\neq i$ with $j\geq2$.
Therefore, if there exists some $i\in \mathbb{N}$ such that $T^*=i T_*$, the kernel space of $D_v F\left(0,T^*\right)$ is just $V_1\cup V_i$.\qed\\

%

We would like to point out that
the kernel space of $D_v F\left(0,T^*\right)$ may not be one-dimensional because $T^*$ may
equal to $jT_*$ for some $j\in \mathbb{N}$ with $j\geq2$.
If $T^*\neq jT_*$ for any $j\in \mathbb{N}$ with $j\geq2$, the kernel space of $D_v F\left(0,T^*\right)$ is also just $V_1$.

\section{Proofs of Theorems 1.1--1.2}

\quad\, We first recall the famous Crandall-Rabinowitz bifurcation theorem \cite[Theorem 1.7 and Theorem 1.18]{Crandall}.
\\ \\
\textbf{Proposition 4.1.} \emph{Let $X$ and $Y$ be real Banach spaces, $U$ be a neighborhood of $0$ in $X$ and $F:\Lambda\times U\longrightarrow Y$ be a $C^1$ map with $F(\lambda,0)=0$ for any $\lambda\in \Lambda$, where $\Lambda$ is an open interval in $\mathbb{R}$. Suppose that}

1. $\text{dim} \text{Ker}\left(D_u F(\mu,0)\right)=\text{codim} \text{Im}\left(D_u F(\mu,0)\right)=1$ \emph{for some $\mu\in \Lambda$},

2. \emph{The partial derivatives $F_\lambda$, $F_u$, $F_{\lambda u}$ exist and are continuous},

3. $D_{\lambda u}F(\mu, 0)\left[w_0\right]\not\in \text{Im}\left(D_u F(\mu,0)\right)$, \emph{where} $w_0\in \text{Ker}\left(D_u F(\mu,0)\right)\setminus\{0\}$.

\noindent\emph{If $X_0$ is any complement of} $\text{Ker}\left(D_u F(\mu,0)\right)$, \emph{then there exist an open interval $I_0 = \left(-\delta_0,\delta_0\right)$ and continuous
functions $\lambda : I_0\rightarrow \mathbb{R}$ and $\psi : I_0 \rightarrow  X_0$ such that $\lambda(0) = \mu$, $\psi(0) = 0$ and $F(\lambda(s), sw_0 + s\psi(s)) = 0$ for
$s\in I_0$ and $F^{-1}\{0\}$ near $(\mu, 0)$ consists precisely of the curves $u =0$ and
$\Gamma = \left\{(\lambda(s), u(s)): s \in I_0\right\}$.}\\

We next show Theorem 1.1 by verifying the hypotheses of Proposition 4.1.
\\ \\
\textbf{Proof of Theorem 1.1.} From Proposition 3.3, we know that the kernel of the
linearized operator $D_vF\left(0,T_*\right)$ is one-dimensional and is spanned by the function $\cos (t)$.
As that of \cite[Proposition 3.2]{Sicbaldi} we can show that $D_vF\left(0,T_*\right)$ is a formally self-adjoint, first order elliptic operator.
It follows that $D_vF\left(0,T_*\right)$ has closed range. Therefore, $D_vF\left(0,T_*\right)$ is a Fredholm operator of index zero (refer to \cite{Kubrusly}). So its codimension is equal to $1$.
In view of Proposition 3.2, we obtain
\begin{equation}
D_{Tv}F\left(0,T_*\right)\cos(t)=\sigma'\left(T_*\right)\cos(t)\not\in \text{Im}\left(D_vF\left(0,T_*\right)\right).\nonumber
\end{equation}
\indent Applying Proposition 3.1 to $F(v,T)=0$, we obtain that there exist an open interval $I= \left(-\varepsilon,\varepsilon\right)$ and continuous
functions $T : I\rightarrow \mathbb{R}$ and $w : I \rightarrow  \text{Im}\left(D_v F\left(0,T_*\right)\right)$ such that $T(0) = T_*$, $w(0) = 0$ and $F(s\cos(t) + sw(s),T(s)) = 0$ for
$s\in I$ and $F^{-1}\{0\}$ near $\left(0,T_*\right)$ consists precisely of the curves $v =0$ and
$\Gamma = \left\{(v(s),T(s)): s \in I\right\}$. Therefore, for each $s\in  (-\varepsilon,\varepsilon)$, (\ref{overdeterminedeigenvalueproblem}) has a $T(s)$-periodic solution
$u \in  \mathcal{C}^{2,\alpha}\left(\Omega_s\right)$ with the expected sign-changing property on the modified cylinder
\begin{equation}
\Omega_s=\left\{(x,t)\in \mathbb{R}^N\times \mathbb{R}:r(x)<1+s\cos \left(\frac{2\pi }{T(s)}t\right)+s w(s)\left(\frac{2\pi }{T(s)}t\right)\right\},\nonumber
\end{equation}
which is the desired conclusion.
\qed\\

In order to prove Theorem 1.2, we establish a bifurcation result from high-dimensional kernel space in the following.
\\
\\
\textbf{Proposition 4.2.} \emph{Let $X$ and $Y$ be real Banach spaces, $U$ be a neighborhood of $0$ in $X$ and $F:\Lambda\times U\longrightarrow Y$ be a $C^1$ map with $F(\lambda,0)=0$ for any $\lambda\in \Lambda$, where $\Lambda$ is an open interval in $\mathbb{R}$. Suppose that}

1. $\text{dim} \text{Ker}\left(D_u F(\mu,0)\right)=k$ \emph{for some $\mu\in \Lambda$ and any $k\in \mathbb{N}$, where the basis vectors of $\text{Ker}\left(D_u F(\mu,0)\right)$ are denoted by $\left\{w_i\right\}_{i=1}^k$};

2. \emph{For any fixed $i\in\{1,\ldots,k\}$, let $X_0$ be any complement of $\text{Ker}\left(D_u F(\mu,0)\right)$ and
$\widetilde{F}$ be the restriction of $F$ on $\Lambda\times \widetilde{U}$ with $\widetilde{U}$ being the restriction of $U$ on $\widetilde{X}:= \text{span}\left\{w_i\right\}\oplus X_0$. Assume that
$\widetilde{Y}\subseteq Y$ is a closed subspace of $Y$ such that the image of $\widetilde{F}$ is contained in $\widetilde{Y}$ and
$\text{dim} \text{Ker}\left(D_u \widetilde{F}(\mu,0)\right)=\text{codim} \text{Im}\left(D_u \widetilde{F}(\mu,0)\right)$};

3. \emph{The double cross derivative $\widetilde{F}_{\lambda u}$ exists and is continuous};

4. $D_{\lambda u}\widetilde{F}(\mu, 0)\left[w_i\right]\not\in \text{Im}\left(D_u \widetilde{F}(\mu,0)\right)$;

\noindent\emph{Then there exist an open interval $I_0 = \left(-\delta_0,\delta_0\right)$ and continuous
functions $\lambda : I_0\rightarrow \mathbb{R}$ and $\psi : I_0 \rightarrow  X_0$ such that $\lambda(0) = \mu$, $\psi(0) = 0$ and $F(\lambda(s), sw_i + s\psi(s)) = 0$ for
$s\in I_0$ and $F^{-1}\{0\}$ near $(\mu, 0)$ consists precisely of the curves $u =0$ and
$\Gamma = \left\{(\lambda(s), u(s)): s \in I_0\right\}$.}\\

When $k=1$, the conclusion of Proposition 4.2 is just the famous Crandall-Rabinowitz bifurcation theorem, where $\widetilde{X}=X$ and $\widetilde{Y}=Y$.
Unlike the case in \cite{Westreich}, here we do not require that $k$ is odd. Moreover, this conclusion is better than that of \cite{Westreich} where the only bifurcation point was obtained.
\\
\\
\textbf{Proof of Proposition 4.2.} We can find that
\begin{equation}
\text{dim} \text{Ker}\left(D_u \widetilde{F}(\mu,0)\right)=1.\nonumber
\end{equation}
Since $D_u \widetilde{F}\left(\lambda,0\right)$ is a Fredholm operator with index zero, there exists a one-dimensional
closed subspace $\widetilde{Z}$ such that
\begin{equation}
\widetilde{Y}=\widetilde{Z}\oplus \text{Im}\left(D_u \widetilde{F}(\mu,0)\right).\nonumber
\end{equation}
For any fixed $i\in\{1,\ldots,k\}$, define $G:\mathbb{R}\times X_0\times \mathbb{R}\longrightarrow \widetilde{Y}$ by
\begin{equation}
G(s,z,\lambda)=\left\{
\begin{array}{ll}
\frac{1}{s}\widetilde{F}\left(\lambda,s\left(w_i+z\right)\right)\,\, &s\neq0,\\
D_u \widetilde{F}\left(\lambda,0\right)\left(w_i+z\right) &s=0.
\end{array}
\right.\nonumber
\end{equation}
Clearly, $G\left(0,0,\mu\right)=0$.
Since $\widetilde{F}$ is $C^1$ and $\widetilde{F}_{\lambda u}$ is continuous, it follows that $G_z$ and $G_\lambda$ are continuous
and
\begin{equation}
G_z\left(0,0,\mu\right)=D_u \widetilde{F}(\mu,0)\,\,\,\text{and}\,\,\,G_\lambda\left(0,0,\lambda_0\right)=\widetilde{F}_{\lambda u}\left(\mu,0\right)w_i.\nonumber
\end{equation}
Then we have that
\begin{equation}
G_{(z, \lambda)}\left(0,0,\mu\right)\left(z_*,\lambda_*\right)=D_u \widetilde{F}(\mu,0)z_*+\lambda_*\widetilde{F}_{\lambda u}\left(\mu,0\right)w_i:X_0\times \mathbb{R}\longrightarrow \widetilde{Y}.\nonumber
\end{equation}
We next prove that $G_{(z, \lambda)}\left(0,0,\mu\right)$ is an isomorphism on to $\widetilde{Y}$.

We first assume that
\begin{equation}
G_{(z, \lambda)}\left(0,0,\mu\right)\left(z_*,\lambda_*\right)=D_u \widetilde{F}(\mu,0)z_*+\lambda_*\widetilde{F}_{\lambda u}\left(\mu,0\right)w_i=0.\nonumber
\end{equation}
If $\lambda_*$ does not vanish, then we have that
\begin{equation}
\widetilde{F}_{\lambda u}\left(\mu,0\right)w_i=-\lambda_*^{-1}D_u \widetilde{F}(\mu,0)z_*\in \text{Im}\left(D_u \widetilde{F}(\mu,0)\right),\nonumber
\end{equation}
which is contradicted with the transversality condition.
Thus we have that
\begin{equation}
\lambda_*=0,\nonumber
\end{equation}
which leads further to that $D_u \widetilde{F}(\mu,0)z_*=0$. Since $D_u \widetilde{F}(\mu,0):X_0\rightarrow \text{Im}\left(D_u \widetilde{F}(\mu,0)\right)$ is an isomorphism, we obtain that $z_*=0$.
Hence we have shown that $G_{(z, \lambda)}\left(0,0,\mu\right)$ is an injection.

We now show that $G_{(z, \lambda)}\left(0,0,\mu\right)$ is also a surjection.
For any $y\in \widetilde{Y}$ such that
\begin{equation}
G_{(z, \lambda)}\left(0,0,\mu\right)\left(z_*,\lambda_*\right)=D_u \widetilde{F}(\mu,0)z_*+\lambda_*\widetilde{F}_{\lambda u}\left(\mu,0\right)w_i=y.\nonumber
\end{equation}
From the above argument we know that $\widetilde{F}_{\lambda u}\left(\mu,0\right)w_i\in \widetilde{Y}$.
By the Hahn-Banach theorem,
there exists a linear functional $l\in \widetilde{Y}^*$ with $\widetilde{Y}^*$ being the dual space of $\widetilde{Y}$ such that
\begin{equation}
l\left(\widetilde{F}_{\lambda u}\left(\mu,0\right)w_i\right)=1,\,\,\, \widetilde{Y}=\text{span}\left\{\widetilde{F}_{\lambda u}\left(\mu,0\right)w_i\right\}\oplus \widetilde{Y}_0, \nonumber
\end{equation}
where $\widetilde{Y}_0=\{v\in \widetilde{Y}:l(v)=0\}=\text{Im}\left(D_u \widetilde{F}(\mu,0)\right)$.
Applying $l$ on the both sides of $G_{(z, \lambda)}\left(0,0,\mu\right)\left(z_*,\lambda_*\right)=y$, we obtain that
$\lambda_*=l(y)$.
Furthermore, we obtain that
\begin{equation}
D_u \widetilde{F}(\mu,0)z_*=y-l(y)\widetilde{F}_{\lambda u}\left(\mu,0\right)w_i.\nonumber
\end{equation}
Since $D_u \widetilde{F}(\mu,0):X_0\rightarrow \text{Im}\left(D_u \widetilde{F}(\mu,0)\right)$ is an isomorphism, it is reversible. Then we have that
\begin{equation}
z_*=\left(D_u \widetilde{F}(\mu,0)\right)^{-1}\left(y-l(y)\widetilde{F}_{\lambda u}\left(\mu,0\right)w_i\right).\nonumber
\end{equation}
Therefore, we prove that $G_{(z, \lambda)}\left(0,0,\mu\right)$ is a surjection.

Based on the argument above, it is obvious that the inverse mapping of $G_{(z, \lambda)}\left(0,0,\mu\right)$ is continuous.
Since
\begin{equation}
G_{(z, \lambda)}\left(0,0,\mu\right):X_0\times \mathbb{R}\longrightarrow \widetilde{Y}\nonumber
\end{equation}
is continuous and an one-to-one mapping, we conclude that it is also a homeomorphic mapping.
Thus, the desired conclusions can be obtained by applying the
Implicit function theorem.\qed
~\\

Further, choosing a combination-type basis in kernel space, we have the following Corollary.
\\ \\
\textbf{Corollary 4.1.} \emph{Let $X$ and $Y$ be real Banach spaces, $U$ be a neighborhood of $0$ in $X$ and $F:\Lambda\times U\longrightarrow Y$ be a $C^1$ map with $F(\lambda,0)=0$ for any $\lambda\in \Lambda$, where $\Lambda$ is an open interval in $\mathbb{R}$. Suppose that}

1. $\text{dim} \text{Ker}\left(D_u F(\mu,0)\right)=k$ \emph{for some $\mu\in \Lambda$ and any $k\in \mathbb{N}$, where the basis vectors of $\text{Ker}\left(D_u F(\mu,0)\right)$ are denoted by $\left\{w_i\right\}_{i=1}^k$};

2. \emph{For any fixed nonzero real array $\left\{t_i\right\}$ with $\sum_{i=1}^kt_i^2=1$, let $X_0$ be any complement of $\text{Ker}\left(D_u F(\mu,0)\right)$ and
$\widetilde{F}$ be the restriction of $F$ on $\Lambda\times \widetilde{U}$ with $\widetilde{U}$ being the restriction of $U$ on $\widetilde{X}:= \text{span}\left\{\sum_{i=1}^kt_iw_i\right\}\oplus X_0$. Assume that
$\widetilde{Y}\subseteq Y$ is a closed subspace of $Y$ such that the image of $\widetilde{F}$ is contained in $\widetilde{Y}$ and
$\text{dim} \text{Ker}\left(D_u \widetilde{F}(\mu,0)\right)=\text{codim} \text{Im}\left(D_u \widetilde{F}(\mu,0)\right)$};

3. \emph{The double cross derivative $\widetilde{F}_{\lambda u}$ exists and is continuous};

4. $D_{\lambda u}\widetilde{F}(\mu, 0)\left[\sum_{i=1}^kt_iw_i\right]\not\in \text{Im}\left(D_u \widetilde{F}(\mu,0)\right)$;

\noindent\emph{Then there exist an open interval $I_0 = \left(-\delta_0,\delta_0\right)$ and continuous
functions $\lambda : I_0\rightarrow \mathbb{R}$ and $\psi : I_0 \rightarrow  X_0$ such that $\lambda(0) = \mu$, $\psi(0) = 0$ and $F(\lambda(s), s\sum_{i=1}^kt_iw_i + s\psi(s)) = 0$ for
$s\in I_0$ and $F^{-1}\{0\}$ near $(\mu, 0)$ consists precisely of the curves $u =0$ and
$\Gamma = \left\{(\lambda(s), u(s)): s \in I_0\right\}$.}
\\ \\
\textbf{Proof.} By some rotational transformation $T$, $\widetilde{X}$ can change into $\widehat{X}=\text{span}\left\{w_i\right\}\oplus T\left(X_0\right)$.
Then, applying Proposition 4.2 on $\widehat{X}$ and using the reverse action of $T$, we can obtained the desired conclusion.\qed\\

We now show Theorem 1.2 by verifying the hypotheses of Corollary 4.1.
\\ \\
\textbf{Proof of Theorem 1.2.} If $T^*\neq mT_*$ for any $m\in \mathbb{N}$, the kernel space of $D_v F\left(0,T^*\right)$ is $V_1$.
Then, in view of $\sigma'\left(T^*\right)>0$, repeating the argument as that of Theorem 1.1 we have the desired conclusion.

If there exists some $m\in \mathbb{N}$ such that $T^*=m T_*$, it follows from Proposition 3.3 that the kernel space of $D_v F\left(0,T^*\right)$ is two-dimensional and is spanned by the functions $\cos (t)$ and $\cos(mt)$ with some $m\geq2$. In addition, it is easy to check that
\begin{equation}
\mathcal{C}^{2,\alpha}_{\text{even},0}\left(\mathbb{R}/2\pi \mathbb{Z}\right)=\text{span}\{\cos(t),\cos(mt)\}\oplus X_0,\nonumber
\end{equation}
where $X_0$ is spanned by $\{\cos(ix)\}_{i}$ with $i\in \mathbb{N}$ and $i\neq 1,m$. For any given nonzero constants $\beta$ and $\gamma$ with $\beta^2+\gamma^2=1$, we define
\begin{equation}
\widetilde{X}=\text{span}\{\beta\cos(t)+\gamma\cos(mt)\}\oplus X_0\nonumber
\end{equation}
and let
\begin{equation}
\widetilde{F}: \mathbb{R^+}\times\widetilde{X}\rightarrow \widetilde{Y}, \nonumber
\end{equation}
where $\widetilde{Y}$ is a closed subspace of $\mathcal{C}^{1,\alpha}_{\text{even},0}\left(\mathbb{R}/2\pi \mathbb{Z}\right)$.
Reasoning as that of Theorem 1.1, $D_v\widetilde{F}\left(0,T^*\right)$ is a formally self-adjoint, first order elliptic operator and has closed range. So $D_v\widetilde{F}\left(0,T^*\right)$ is also a Fredholm operator of index zero (refer to \cite{Kubrusly}) with its codimension is equal to $1$.

In addition, it is easy to check that the double cross derivative $\widetilde{F}_{\lambda u}$ exists and is continuous. At last, let us verify the transversality condition.
By using the relation $\sigma_m\left(T\right)=\sigma\left(T/m\right)$, we have that
\begin{equation}
D_{Tv}\widetilde{F}\left(0,T^*\right)\left(\beta\cos(t)+\gamma\cos(mt)\right)=\beta\sigma'\left(T^*\right)\cos(t)+\frac{\gamma}{m}\sigma'\left(T_*\right)\cos(mt),\nonumber
\end{equation}
where $\sigma'\left(T^*\right)>0$ and $\sigma'\left(T_*\right)>0$.
For any $v$ belonging to $\widetilde{X}$ with $v>-1$, using the Fourier expansion, $v$ can be written as
\begin{equation}
v=\sum_{i\geq1}a_i\cos(it),\nonumber
\end{equation}
where $a_1=\beta$ and $a_m=\gamma$.
We know that
\begin{equation}
D_v\widetilde{F}\left(0,T\right)v=\sum_{i\geq1}\sigma_i\left(T\right)a_i\cos(it)=\sum_{i\geq1}\sigma\left(\frac{T}{i}\right)a_i\cos(it).\nonumber
\end{equation}
Thus, in view of $\sigma\left(T^*\right)=0$ and $\sigma_m\left(T^*\right)=\sigma\left(T_*\right)=0$, we have
\begin{equation}
D_v\widetilde{F}\left(0,T^*\right)v=\sum_{i\geq2,i\neq m}\sigma_i\left(T^*\right)a_i\cos(it)=\sum_{i\geq2,i\neq m}\sigma\left(\frac{T^*}{i}\right)a_i\cos(it).\nonumber
\end{equation}
Using Proposition 3.2, we deduce that
\begin{equation}
\sigma\left(\frac{T^*}{i}\right)\neq0\nonumber
\end{equation}
for any $i\geq2$, $i\neq m$.
Consequently, the image of $D_v\widetilde{F}\left(0,T^*\right)$ is the closure of
\begin{equation}
\bigoplus_{i\geq2, i\neq m}V_i\nonumber
\end{equation}
in $\widetilde{Y}$.
Then it follows the fact that
\begin{equation}
\beta\sigma'\left(T^*\right)\cos(t)+\frac{\gamma}{m}\sigma'\left(T_*\right)\cos(mt) \not \in\text{Im}\left(D_v\widetilde{F}\left(0,T^*\right)\right),\nonumber
\end{equation}
which means that
\begin{equation}
D_{Tv}\widetilde{F}\left(0,T^*\right)\left(\beta\cos(t)+\gamma\cos(mt)\right) \not \in\text{Im}\left(D_v\widetilde{F}\left(0,T^*\right)\right).\nonumber
\end{equation}

By applying Corollary 4.1 to $\widetilde{F}(v,T)=0$, we obtain that there exist an open interval $I= \left(-\varepsilon,\varepsilon\right)$ and continuous
functions $T : I\rightarrow \mathbb{R}$, $w : I \rightarrow  \text{Im}\left(D_v F\left(0,T^*\right)\right)$ such that $T(0) = T^*$, $w(0) = 0$ and $F(s(\beta\cos(t)+\gamma\cos(mt)) + sw(s),T(s)) = 0$ for
$s\in I$ with $\beta^2+\gamma^2=1$ and $F^{-1}\{0\}$ near $\left(0,T^*\right)$ consists precisely of the curves $v =0$ and
$\Gamma = \left\{(v(s),T(s)): s \in I\right\}$. Therefore, for each $s\in  (-\varepsilon,\varepsilon)$, (\ref{overdeterminedeigenvalueproblem}) has a $T(s)$-periodic solution
$u \in  \mathcal{C}^{2,\alpha}\left(\Omega_s\right)$ with the expected sign-changing property on the modified cylinder
\begin{equation}
\Omega_s=\left\{(x,t)\in \mathbb{R}^{N+1}:\vert x\vert<1+s\left(\beta\cos \left(\frac{2\pi }{T(s)}t\right)+\gamma\cos \left(\frac{2m\pi }{T(s)}t\right)\right)+s w(s)\left(\frac{2\pi }{T(s)}t\right)\right\},\nonumber
\end{equation}
as desired.\qed

\section{The case of $N=1$}

\quad\, Note that the arguments of Propositions 3.1--3.2 may not be valid for the one-dimensional case.
For $N=1$, it is necessary to consider the zeros of $J_{-\frac{1}{2}}(s)$ or $J_{-\frac{3}{2}}(s)$.
While, we cannot study the zero distribution of $J_{-\frac{1}{2}}(s)$ or $J_{-\frac{3}{2}}(s)$ with the interlace property or eigenfunction because $-3/2<-1/2<0$. Moreover, from the argument of Proposition 3.2, we find that $\sigma(\mu)<0$ holds only under the condition $N\geq2$.
In fact, we will see that $\sigma(\mu)=0$ if $N=1$.
Thus, we need new methods to study the one-dimensional case. Fortunately, in the one-dimensional case, many conclusions can be calculated specifically.

For $N=1$, we have that
\begin{equation}
\lambda_1=\frac{\pi^2}{4},\,\,\lambda_2=\frac{9\pi^2}{4}\nonumber
\end{equation}
and
\begin{equation}
\phi_2=\frac{1}{\sqrt{2\pi}}\cos\left(\frac{3\pi}{2}r\right).\nonumber
\end{equation}
So we have that
\begin{equation}
\phi_2'(1)=\frac{3\sqrt{2\pi}}{4}\,\,\,\text{and}\,\,\,\phi_2''(1)=0.\nonumber
\end{equation}
It follows that
\begin{equation}
\sigma(T)=c'(1)+\phi_2''(1)=c'(1),\nonumber
\end{equation}
$c:=c_1$ is the continuous solution on $[0,1]$ of
\begin{equation}\label{n=1ckequation}
\left(\partial_r^2+\lambda_2\right)c-\left(\frac{2\pi}{T}\right)^2c=0
\end{equation}
with $c(1)=-\phi_2'(1)$ and $c'(0)=0$. Proposition 2.1 implies $\sigma(T)$ is analytic when $T\neq\sqrt{2}$.
\\ \\
\textbf{Lemma 5.1.} \emph{The function $\sigma:\left(0,\sqrt{2}\right)\cup \left(\sqrt{2},+\infty\right)\rightarrow \mathbb{R}$ has exactly
two zeros $4/3$ and $4\sqrt{5}/5$ such that
$\sigma'\left(4/3\right)>0$ and $\sigma'\left(4\sqrt{5}/5\right)>0$. Moreover, $\sigma(T)$ satisfies
$$\lim_{T\rightarrow 0^+}\sigma(T)=-\infty, \quad  \lim_{T\rightarrow \left(\sqrt{2}\right)^-}\sigma(T)=+\infty$$
and
$$\lim_{T\rightarrow \left(\sqrt{2}\right)^+}\sigma(T)=-\infty, \quad \lim_{T\rightarrow+\infty}\sigma(T)=+\infty.$$}
\\ \\
\textbf{Proof.} Let
\begin{equation}
\alpha(T)=\frac{9\pi^2}{4}-\left(\frac{2\pi}{T}\right)^2.\nonumber
\end{equation}
Then we have that
\begin{equation}
\alpha\left(\frac{4}{3}\right)=0, \,\,\,\alpha\left(\sqrt{2}\right)=\frac{\pi^2}{4}.\nonumber
\end{equation}
The solution of (\ref{n=1ckequation}) is
\begin{equation}
c(r)=\left\{
\begin{array}{ll}
-\frac{3\sqrt{2\pi}}{4}\frac{\cosh\left(\sqrt{-\alpha(T)}r\right)}{\cosh\left(\sqrt{-\alpha(T)}\right)}\,\, &\text{if}\,\, T\in\left(0,\frac{4}{3}\right),\\
-\frac{3\sqrt{2\pi}}{4}&\text{if}\,\, T=\frac{4}{3},\\
-\frac{3\sqrt{2\pi}}{4}\frac{\cos\left(\sqrt{\alpha(T)}r\right)}{\cos\left(\sqrt{\alpha(T)}\right)} &\text{if}\,\, T\in\left(\frac{4}{3},\sqrt{2}\right)\cup\left(\sqrt{2},+\infty\right).
\end{array}
\right.\nonumber
\end{equation}
Then we get that
\begin{equation}
\sigma(T)=c'(1)=\left\{
\begin{array}{ll}
-\frac{3\sqrt{2\pi}}{4}\sqrt{-\alpha(T)}\tanh\left(\sqrt{-\alpha(T)}\right)\,\, &\text{if}\,\, T\in\left(0,\frac{4}{3}\right),\\
0&\text{if}\,\, T=\frac{4}{3},\\
\frac{3\sqrt{2\pi}}{4}\sqrt{\alpha(T)}\tan\left(\sqrt{\alpha(T)}\right) &\text{if}\,\, T\in\left(\frac{4}{3},\sqrt{2}\right)\cup\left(\sqrt{2},+\infty\right).
\end{array}
\right.\nonumber
\end{equation}
It follows that
\begin{equation}
\sigma(T)<0 \,\,\,\text{for}\,\,\,T<\frac{4}{3}\nonumber
\end{equation}
and
\begin{equation}
\sigma(T)>0 \,\,\,\text{for}\,\,\,T\in\left(\frac{4}{3},\sqrt{2}\right).\nonumber
\end{equation}
The unique zero of $\sigma$ in $\left(0,\sqrt{2}\right)$ is $4/3$.
It is obvious that $\lim_{x\rightarrow \left(\frac{\pi}{2}\right)^-}\tan(x)=+\infty$ and $\lim_{x\rightarrow \left(\frac{\pi}{2}\right)^+}\tan(x)=-\infty$, we get that
$\lim_{T\rightarrow \left(\sqrt{2}\right)^-}\sigma(T)=+\infty$ and $\lim_{T\rightarrow \left(\sqrt{2}\right)^+}\sigma(T)=-\infty$.
Since $\alpha(T)$ goes to $-\infty$ as $T\rightarrow0$, one can see that $\lim_{T\rightarrow 0^+}\sigma(T)=-\infty$.
Further, in view of $\lim_{x\rightarrow \left(\frac{3\pi}{2}\right)^+}\tan x=+\infty$ and $\lim_{T\rightarrow+\infty}\alpha(T)=9\pi^2/4$, we have that
$\lim_{T\rightarrow+\infty}\sigma(T)=+\infty$.
Therefore, $\sigma(T)$ also has at least one zero in $\left(\sqrt{2},+\infty\right)$.

We compute that
\begin{equation}
\sigma'(T)=\left\{
\begin{array}{ll}
\frac{3\sqrt{2\pi}\alpha'(T)}{8\sqrt{-\alpha(T)}}\left(\tanh\left(\sqrt{-\alpha(T)}\right)+\frac{\sqrt{-\alpha(T)}}{\cosh^2\left(\sqrt{-\alpha(T)}\right)}\right)\,\, &\text{if}\,\, T\in\left(0,\frac{4}{3}\right),\\
\frac{3\sqrt{2\pi}\alpha'(T)}{8\sqrt{\alpha(T)}}\left(\tan\left(\sqrt{\alpha(T)}\right)+\sqrt{\alpha(T)}\sec^2\left(\sqrt{\alpha(T)}\right)\right) &\text{if}\,\, T\in\left(\frac{4}{3},\sqrt{2}\right)\cup\left(\sqrt{2},+\infty\right).
\end{array}
\right.\nonumber
\end{equation}
It follows that $\sigma'(T)>0$ for $T\in\left(0,4/3\right)$. We next compute the derivative of $\sigma$ at $T=4/3$. 

The left derivative of $\sigma$ at $T=4/3$ is
\begin{eqnarray}
\sigma'\left(\frac{4}{3}-0\right)&=&\lim_{T\rightarrow\left(\frac{4}{3}\right)^-}\frac{3\sqrt{2\pi}\alpha'(T)}{8\sqrt{-\alpha(T)}}
\left(\tanh\left(\sqrt{-\alpha(T)}\right)+\frac{\sqrt{-\alpha(T)}}{\cosh^2\left(\sqrt{-\alpha(T)}\right)}\right)\nonumber\\
&=&\frac{81\pi^2\sqrt{2\pi}}{64}\lim_{T\rightarrow\left(\frac{4}{3}\right)^-}
\left(\frac{\tanh\left(\sqrt{-\alpha(T)}\right)}{\sqrt{-\alpha(T)}}+\frac{1}{\cosh^2\left(\sqrt{-\alpha(T)}\right)}\right),\nonumber
\end{eqnarray}
where we have used the fact of $\alpha'\left(4/3\right)=27\pi^2/8$.
By the L'Hospital rule we know that
\begin{equation}
\lim_{x\rightarrow0}\frac{\tanh(x)}{x}=\lim_{x\rightarrow0}\frac{1}{\cosh^2(x)}=1.\nonumber
\end{equation}
Thus, we have that
\begin{equation}
\lim_{T\rightarrow\left(\frac{4}{3}\right)^-}
\left(\frac{\tanh\left(\sqrt{-\alpha(T)}\right)}{\sqrt{-\alpha(T)}}+\frac{1}{\cosh^2\left(\sqrt{-\alpha(T)}\right)}\right)=2\nonumber
\end{equation}
due to $\lim_{T\rightarrow\left(\frac{4}{3}\right)^-}\alpha(T)=0$.
Therefore, the left derivative of $\sigma$ at $T=4/3$ is
\begin{eqnarray}
\sigma'\left(\frac{4}{3}-0\right)=\frac{81\sqrt{2}\pi^{\frac{5}{2}}}{32}.\nonumber
\end{eqnarray}
Similarly, the right derivative of $\sigma$ at $T=4/3$ is
\begin{eqnarray}
\sigma'\left(\frac{4}{3}-0\right)&=&\lim_{T\rightarrow\left(\frac{4}{3}\right)^-}
\frac{3\sqrt{2\pi}\alpha'(T)}{8\sqrt{\alpha(T)}}\left(\tan\left(\sqrt{\alpha(T)}\right)+\sqrt{\alpha(T)}\sec^2\left(\sqrt{\alpha(T)}\right)\right)\nonumber\\
&=&\frac{81\pi^2\sqrt{2\pi}}{64}\lim_{T\rightarrow\left(\frac{4}{3}\right)^-}
\left(\frac{\tan\left(\sqrt{\alpha(T)}\right)}{\sqrt{\alpha(T)}}+\frac{1}{\cos^2\left(\sqrt{\alpha(T)}\right)}\right),\nonumber\\
&=&\frac{81\sqrt{2}\pi^{\frac{5}{2}}}{32}.\nonumber
\end{eqnarray}
So $\sigma$ is derivable at $T=4/3$ and $\sigma'\left(4/3\right)>0$.

We now consider the case $T\in\left(4/3,\sqrt{2}\right)$. In this case, we see that $\alpha\in\left(0,\pi^2/4\right)$. It implies that $\sigma'(T)>0$. We finally prove that
$\sigma'(T)>0$ for $T\in\left(\sqrt{2},+\infty\right)$. When $T>\sqrt{2}$, we see that $\sqrt{\alpha}\in\left(\pi/2,3\pi/2\right)$. For $x\in\left(\pi/2,3\pi/2\right)$, let
\begin{eqnarray}
f(x)=\tan(x)+x\sec^2(x).\nonumber
\end{eqnarray}
Then we find that
\begin{eqnarray}
f(x)=\frac{\sin(2x)+2x}{2\cos^2(x)}.\nonumber
\end{eqnarray}
Let $g(x)=\sin(2x)+2x$. We get that $g(\pi/2)=\pi$ and $g'(x)=2(\cos(2x)+1)>0$ for $x\in\left(\pi/2,3\pi/2\right)$.
So, $g(x)>0$, and it further implies $f(x)>0$. Hence, we get that $\sigma'(T)>0$ for $T\in\left(\sqrt{2},+\infty\right)$.

Thus, we obtain that  $\sigma'(T)>0$ for $T\in\left(0,\sqrt{2}\right)\cup\left(\sqrt{2},+\infty\right)$.
Hence, $\sigma$ also has a unique zero in $\left(\sqrt{2},+\infty\right)$ which is denoted by $T^*$ such that $\sigma'\left(T^*\right)>0$.
In fact, from the expression of $\sigma(T)$ we derive that $T^*=4\sqrt{5}/5$, which is not a multiple of $T_*=4/3$ in this case.
\qed\\

In view of Lemma 5.1, it is easy to verify the conditions in Propositions 4.1--4.2. Thus the conclusions of Theorems 1.1--1.2 are all valid for $N=1$.

\section{Appendix}

In \cite{Schlenk}, Schlenk and Sicbaldi gave the following claim.
\\ \\
\textbf{Claim A.1.} \emph{One has $J_\nu^2(s)>J_{\nu-1}(s)J_{\nu+1}(s)$ for all $s\in\left(0,j_{\nu,1}\right)$.}
\\

This claim plays a key role to show $\sigma'(T)\neq0$ which is the most important step in verifying the transversality condition.
Unfortunately, it contains a small gap in the argument.
They used the following relations
\begin{equation}\label{wrongrelations}
j_{\nu-1,1}<j_{\nu,1}<j_{\nu+1,1},\,\,\,j_{\nu,1}<j_{\nu-1,2},
\end{equation}
which is just interlace property of zeros.
However, when $\nu-1<0$, this property may no longer hold.
For example, in the case of $N=2$, we have that $\nu-1=-1$ and
$J_{-1}(s)=-J_1(s)$.
By interlace property of zeros, we know that $j_{0,1}<j_{1,1}<j_{0,2}<j_{1,2}$.
So, $j_{-1,1}=j_{1,1}\in\left(j_{0,1},j_{0,2}\right)$ and $j_{-1,2}=j_{1,2}>j_{0,2}$, which indicate that relations (\ref{wrongrelations}) are not correct.
Here we reinvestigate this claim by filling the above gap.
\\ \\
\textbf{Proof Claim A.1.} When $N\geq4$, we see that $\nu-1\geq0$, the interlace property holds. Therefore, the argument of \cite{Schlenk} is valid for $N\geq4$ and it remains to prove the case of $N=2,3$.

We first consider the case of $N=2$.
In this case, $\nu=0$ and $J_{-1}(s)=-J_1(s)$.
It follows that $J_{\nu-1}(s)J_{\nu+1}(s)=-J_1^2(s)\leq0$.
Hence we have that $J_{\nu}^2(s)-J_{\nu-1}(s)J_{\nu+1}(s)>0$ in $\left(0,j_{\nu,1}\right)$.

We next assume $N=3$. In this case, by Lemma 2.1, $j_{\nu,1}=\pi$ and $j_{\nu,2}=2\pi$ with $\nu=1/2$.
We have known that
$J_{-1/2}(s)$ is positive in $\left(0,\pi/2\right)$ and the first positive zero is $\pi/2$.
By the interlace property of zeros, we have that $j_{\nu+1,1}>\pi$.
Hence $J_{-1/2}(s)$ and $J_{3/2}(s)$ have
opposite sign in $\left(\pi/2,\pi\right)$.
It follows that $J_{\nu}^2(s)-J_{\nu-1}(s)J_{\nu+1}(s)>0$ in $\left(\pi/2,\pi\right)$.
At $s=\pi/2$, $J_{\nu+1}(s)>0$, $J_{\nu}(s)>0$ and $J_{\nu-1}(s)=0$. Thus, $J_{\nu}^2(\pi/2)-J_{\nu-1}(\pi/2)J_{\nu+1}(\pi/2)>0$.

It suffices to study the case of $s\in\left(0,\pi/2\right)$.
In this interval, $J_{-1/2}(s)$ and $J_{3/2}(s)$ are all positive.
Reasoning as that of Proposition 3.2 we can show that $J_{\nu}^2(s)-J_{\nu-1}(s)J_{\nu+1}(s)$ is strictly increasing in $\left(0,\pi/2\right)$.
We have shown that $J_{\nu}^2(s)-J_{\nu-1}(s)J_{\nu+1}(s)>0$ at the right endpoint $\pi/2$.
By the asymptotic formula (\ref{relationsforbesselfunctionz=0}), we have that
\begin{equation}
\lim_{s\searrow0}\left(J_{\nu}^2(s)-J_{\nu-1}(s)J_{\nu+1}(s)\right)=0.\nonumber
\end{equation}
Therefore, we can obtain $J_{\nu}^2(s)-J_{\nu-1}(s)J_{\nu+1}(s)>0$ in $\left(0,j_{\nu,1}\right)$.
\qed
\\ \\
\textbf{Acknowledgment}
\bigskip\\
\indent The authors would like to express their gratitude to the anonymous referee for his/her careful
reading of the work and many valuable comments and suggestions. The authors also would like to express their gratitude to Professor Tobias Colding for his seriously and responsibly handling the manuscript.

\bibliographystyle{amsplain}
\makeatletter
\def\@biblabel#1{#1.~}
\makeatother


\providecommand{\bysame}{\leavevmode\hbox to3em{\hrulefill}\thinspace}
\providecommand{\MR}{\relax\ifhmode\unskip\space\fi MR }
\providecommand{\MRhref}[2]{%
  \href{http://www.ams.org/mathscinet-getitem?mr=#1}{#2}
}
\providecommand{\href}[2]{#2}

\end{document}